\documentclass[reqno,a4paper]{article}
\usepackage{graphicx}
\usepackage{color}
\usepackage{psfrag}

\usepackage{amsmath,amsfonts}
\usepackage{a4wide}
\usepackage{xspace}
\usepackage{mathrsfs}

\newcommand{\SPDE}{{SPDE}\xspace}
\newcommand{\SPDEs}{{SPDEs}\xspace}
\newcommand{\PDE}{{PDE}\xspace}
\newcommand{\PDEs}{{PDEs}\xspace}
\newcommand{\PDAE}{{PDAE}\xspace}

\newcommand{\SPDAE}{{SPDAE}\xspace}

\newcommand{\Ito}{{It\^o}\xspace}

\newcommand{\LF}{{\Lambda_{\text{fit}}}}

\newcommand{\real}{{\mathbb{R}}}         
\newcommand{\integers}{{\mathbb{Z}}}         
\newcommand{\Dt}{{\Delta t}}         
\newcommand{\Dx}{{\Delta x}}         

\newcommand{\setR}{{\mathbb{R}}}

\newcommand{\uc}{u_{\mathrm{det}}}
\newcommand{\lc}{\lambda_{\mathrm{det}}}
\newcommand{\uh}{\hat u}
\newcommand{\kh}{{\hat k}}
\newcommand{\Dc}{D_C}

\newcommand{\e}[1]{{\mathrm e^{#1}}}

\newcommand{\E}{\mathbb{E}}
\def\expect#1{\mathbb{E}\left(#1\right)}              
\newcommand{\LAM}{\Lambda_{\min}}
\newcommand{\Lg}{\Lambda_\gamma}
\newcommand{\LAMF}{\Lambda_{\min}^{fix}}
\newcommand{\LgF}{\Lambda_\gamma^{fix}}

\newcommand{\figref}[1]{{Figure~\ref{#1}}}
\newcommand{\tabref}[1]{{Table~\ref{#1}}}
\newcommand{\secref}[1]{{Section~\ref{#1}}}

\title{Freezing Stochastic Travelling Waves}

\author{G. J. Lord\thanks{Department of Mathematics and Maxwell
    Institute, Heriot-Watt University, Edinburgh, 
    EH4 1ER, UK, \texttt{g.j.lord@hw.ac.uk}}
\and V. Th{\"u}mmler\thanks{Universit{\"a}t Bielefeld,Bielefeld, Germany,
\texttt{thuemmle@math.uni-bielefeld.de}}
}

\date{DRAFT \today}

\begin{document}
\maketitle

\begin{abstract}

The aim of this paper is to investigate new numerical methods to
compute travelling wave solutions and new ways to estimate
characteristic properties such as wave speed for stochastically forced
partial differential equations. As a
particular example we consider the Nagumo equation with multiplicative
noise which we mainly consider in the Stratonovich sense.
A standard approach to determine the position and hence speed of a
wave is to compute the evolution of a level set. We compare this
approach against an alternative where the 
wave position is found by minimizing the $L^2$ norm against a fixed
profile. This approach can also be used to stop (or freeze) 
the wave and obtain a stochastic partial differential algebraic
equation that we then discretize and solve. Although attractive as it
leads to a smaller domain size it can be numerically unstable due to
large convection terms. 
We compare numerically the different approaches for estimating the wave speed. 
Minimization against a fixed profile works well provided the support of
the reference function is not too narrow.
We then use these techniques to investigate the effect of both \Ito
and Stratonovich noise on the Nagumo equation as correlation length
and noise intensity increases.


\end{abstract}


\section{Introduction}
\label{sec:intro}
The effects of stochastic forcing on solutions of stochastic partial
differential equations (\SPDEs) have recently received a great deal of
attention in applications ranging from material science, atmosphere
modelling to neural science. Of particular interest are the effects of
noise on travelling waves and fronts as these are often physically
important solutions. We use the term travelling wave to include 
fronts and waves and develop in this paper numerical methods to solve
for stochastic versions of these objects. 

We consider the stochastic PDE
\begin{equation}\label{eq:spde0}
   du = \Big[u_{xx} + f(u)\Big] \,dt + g(u,t)\circ dW(t), \qquad
\mbox{given } \quad u(0)=u^0, \qquad  x \in \setR.
\end{equation}
Although the noise term is written here as a  Stratonovich integral we
consider below the noise in both a Stratonovich and \Ito sense.
For ease of exposition we take $u:\setR \times \setR_+ \to \setR$
although similar numerical proceedures have been applied to
systems of PDEs and waves in $\setR^2$, see for example \cite{bst07}.
For additive noise the function $g$ is taken as a constant whereas for
multiplicative noise $g$ depends on $u$. In the case of no noise
($g=0$) we recover the deterministic \PDE 
\begin{equation}\label{eq:pde}
  u_t = u_{xx} + f(u), \qquad \mbox{given } \quad u(0)=u^0, \qquad x\in\setR.
\end{equation}
In the deterministic case the analysis of travelling waves both
analytically and by numerics is a mature field. This is not the case
for \SPDEs where much of the analysis is performed 
for specific equations or for the case of small noise. Indeed with
stochastic forcing existence for all time of these waves is not
guaranteed and the definition of quantities such as wave speed vary
from system to system. 
Typically the position of a stochastic travelling wave 
is determined from the position of a level set. For small
noise the centres of these fronts can be shown to follow a
rescaled Brownian motion, see \cite{Shardlow:00,Brsssco.etal,Funaki95}.
In the case of multiplicative noise the front may exist
for all times and the wave front may have compact support
\cite{Shga94} and 
is well--defined over some time varying interval in space
$[a(t),b(t)]$ and takes stationary values out side of this
interval.  
There is a well developed literature for the stochastic
Fisher--Kolmogorov-Piscounov equation
\cite{Drng+etal:03,Drng.etal:05,Mllr_Swrs:95,Trbe96} with waves
defined in this form. 
Multiplicative noise is seen to change the wave speed of the wave and the
position is seen to diffuse from the mean (or Goldstein mode), for
reviews see \cite{GrciaOjlvo+Sncho,panja04}. 
However, it is not our aim to replicate these results here. 
We investigate different ways to measure 
the wave speed, using the level set approach as well as a new approach
of minimizing the $L^2$ norm of the wave against a fixed profile.
Furthermore we apply different computational techniques to compute time
dependent waves - this inlcudes freezing the wave to stop it from
travelling. 
We extend a numerical method introduced in  \cite{bt04} for
deterministic \PDEs of the form \eqref{eq:pde}.
This method freezes the wave in the computational domain by adding
a convection term to the equation to compensate for the movement of
the wave.
The convection term that gives the speed of the wave is determined
from an extra algebraic condition from the $L^2$ minimization and the wave
speed is explicitly solved for as a time dependent quantity.
Convergence of this method for \PDEs was considered in \cite{thu05} and
stability of the wave considered in \cite{thu06}. 
We use these techniques to obtain new computational results on the
effect of  multiplicative noise on the wave speed of the Nagumo
equation with Stratonovich and \Ito noise. In contrast to
\cite{GrciaOjlvo+Sncho,Armro_etal,PhysRevE.58.5494}, we present new 
numerical results on the effect of spatial correlation length on the
width of the wave and compare different measures of the wave speed. 
These results incude our new idea of comparing the wave to a reference
function. Furthermore we present new results for \Ito and additive
noise.

Numerically we solve for the wave profile and a time
dependent wave speed for \eqref{eq:spde0} which is, in the case of
stochastic forcing, a random variable.
As a specific example to illustrate the computational
method and to compare against existing techniques we consider 
the scalar Nagumo equation \cite{nay62}
\begin{equation} \label{eq:snagumo}
du = \left[u_{x x} + u(1-u)(u-\alpha)\right]dt +  (\nu + \mu u(1-u))
\circ dW.
\end{equation}
With $\nu\neq0$ and $\mu=0$ we have additive noise and $\mu\neq0$ the noise is
multiplicative. 
For multiplicative noise we have that $u=0$ and $u=1$
are stationary and numerical simulations suggest a wave exists between
them \cite{Armro_etal,PhysRevE.58.5494}. The deterministic equation 
\begin{equation}\label{eq:nagumo}
  u_t = u_{x x} + u(1-u)(u-\alpha),\quad u(x,t)\in\setR,\ x\in\setR,\  t > 0,
\end{equation}
is often used for testing algorithms since  travelling wave 
solutions $u(x,t)=\uc(x-c t)$ connecting the stationary points
$u_-=0$, $u_+=1$ of this equation are explicitly known besides other
explicit solutions, such as pulses, sources and sinks
\cite{air85,cg92}.
These travelling wave solutions depend on the nonlinearity and the
leading profile of initial data $u^0$. Define the function $u_k(x)$ by
\begin{equation}
u_{k}(x) =\left(1+ \e{-kx}\right)^{-1}.
\label{eq:uk}
\end{equation}
We use this function to specify both initial data $u^0$ and
reference functions that have different profiles (by varying $k$).
For $\alpha\in (0,1/2]$ there is a unique asymptotic travelling wave
where as for $\alpha\in (-1,0]$ the 
asymptotic profile and the wave speed depends on the leading profile
$e^{k_0 x}$ of the initial data as $x \to \infty$. 
We summarize results below for the deterministic Nagumo equation,
these are found, for example, in  \cite{GrciaOjlvo+Sncho}. 
%
\begin{itemize}
\item For $\alpha \in (0,1/2]$ the solution $u=u_k$, with $k=1/\sqrt{2}$ is 
  asymptotically stable and all initial front data $u^0$ is attracted to
  this wave. The asymptotic wave speed is given 
  by $c = -\sqrt 2\ (\tfrac 1 2 - \alpha)$.
\item For $\alpha \in (-1/2,0]$ if the initial data $u^0=u_{k_0}$ has 
  $k_0\geq  k_*=-\alpha\sqrt{2}$ then the asymptotic 
  speed is given by $c = -\sqrt 2\ (\tfrac 1 2 - \alpha)$.
If $k_0 < k_*$ then the asymptotic wave speed is $\geq  (k_0^2-\alpha)/k_0$. 
\item For $\alpha \in (-1,-1/2]$ if the initial data $u^0=u_{k_0}$ is such
  that $k_0\geq  k_\dag=\sqrt{|\alpha|}$ 
then the asymptotic speed is given by $2k_{\dag}$.
If $k_0 < k_\dag$ then the asymptotic speed is $\geq2k_{\dag}$.
\end{itemize}

The outline of the rest of the paper is as follows.
We review in \secref{sec:2} the computation of travelling waves in the
deterministic case and extend to the case of stochastic forcing. 
We discuss measures of wave
speeds of a stochastic travelling wave and discuss the numerical approximation.
In \secref{sec:results} we
illustrate the numerical method on the Nagumo equation with
multiplicative noise.
We compare solving the stochastic partial differential equation
(SPDE) and the stochastic partial differential
algebraic equation (SPDAE) where the wave is frozen by minimizing the
$L^2$ distance between a reference function and the travelling
wave. We compare the different measures of wave speeds. We illustrate
from numerics that although a feasible method it can lead to numerical
instability.  We investigate the effect of the choice of reference
function $\uh$ in \secref{sec:uref}.
In \secref{sec:ItoStrat} we present new numerical results for \Ito and
Stratonovich multiplicative noise on how the wave speed changes with
noise intensity. We alos present new results on the effect of the
spatial correlation length. Additive noise is considered in
\secref{sec:add} where we again examine the wave speed with noise
intensity and illustrate the SPDAE approach when new travelling waves are
nucleated. 
Finally we consider weaker versions of the stochastic travelling wave
fixed in the computational domain by mean wave speeds and then discuss
the results and computational method. 

\section{Stochastic travelling waves}
\label{sec:2}


In this section we introduce the (stochastic) differential algebraic
equations that we use to define the travelling wave problem. 
We start by reviewing the more familiar deterministic case before
considering the case with stochastic forcing. In both cases we 
reduce the infinite problem to finite dimensions by truncating the
computational domain and discretizing in space. 

\subsection{Deterministic \PDE and discretization}
Let us assume that equation \eqref{eq:pde} has a travelling wave solution
$u$, so that $u$ can be written as 
\begin{equation}\label{eq:trav.wav.1}
  u(x,t) = \uc(\xi),\quad \xi = x- \lc t,\
\end{equation}
where $\uc \in \mathcal C^2_b(\setR, \setR^m)$ denotes the waveform
and $\lc$ its wave speed. 
In a comoving frame $v(\xi,t)=u(\xi - \lc t,t)$ equation
\eqref{eq:pde} reads
\begin{equation}\label{eq:moving.frame}
  v_t =  v_{\xi \xi} + \lc v_\xi + f(v), \quad \xi \in \setR,\quad t \geq 0
\end{equation}
of which the travelling wave $\uc$ is a stationary solution.
Since the wave speed $\lc$ is generally unknown we transform equation
\eqref{eq:pde} into a co-moving frame with unknown position $\gamma(t)$, i.e. we
insert the ansatz  
$v(x,t) = u(x - \gamma(t), t)$ into \eqref{eq:pde}. Then we obtain
\begin{equation}\label{eq:pde.trafo}
  v_t= v_{xx}+ \lambda  v_x +f(v),
\end{equation}
where $\lambda(t) = \gamma'(t)$.
In order to compensate for the additional variable $\lambda$ we add
a so called phase condition 
\begin{equation}\label{eq:pc}
  0=\psi(v, \lambda)
\end{equation}
which together with \eqref{eq:pde.trafo} forms a partial differential
algebraic equation (\PDAE) \cite{bt04}. 
The position $\gamma$ of the wave can then be calculated by integrating
$  \gamma' = \lambda$,$\gamma(0) = 0$ to get
\begin{equation}\label{eq:gamma}
  \gamma(t) = \int_0^t \lambda(s) ds. 
\end{equation}
For the numerical implementation we need to truncate the spatial domain 
from $x\in \setR$ to $x\in[0,L]$ and impose appropriate boundary conditions
such as  Neumann, Dirichlet or projection boundary conditions \cite{thu05}.
We then solve \eqref{eq:pde.trafo} and \eqref{eq:pc} for $x\in
[0,L]$. 
In contrast to traditional deterministic travelling wave computations
where the steady states of \eqref{eq:moving.frame} are solved for 
with appropriate boundary conditions
this method does not rely on $\lambda$ being a constant wave speed.

Thus far we have not discussed the choice of the phase fixing function
$\psi$ in \eqref{eq:pc}. Since the phase condition only selects one
representative out of the infinite family of solutions, there is some
freedom of choice here. 
The simplest phase condition is to align the solution
with respect to a given reference function $\uh$. 
It is natural to take the $L^2$ norm since we then the minimum can be
found by differentiation and equating to zero (see eg \cite{bt07}). 
For two functions $v$ and $w$ to minimize over shifts in space $y$ the
$L^2$ norm: 
$\min_{y} \|v(x,\cdot)-w(x-y,\cdot)\|_2$. 
Differentiating and equating to zero we find that 
$$\int (v(x,\cdot)-w(x-y,\cdot))w_x dx = 0.$$ 
So for the PDE we take the phase condition to be : 
$$
\psi_{\mathrm{fix}}(u) = \langle \uh_x, u-\uh \rangle.
$$
This choice was termed the template fitting method in \cite{rkml03}.

In our numerical simulations we will discretize in space 
using standard uniformly spaced finite 
differences, so we discretize on a finite grid $x_0,...,x_M$, $u(x_j) = u_j$. 
For the second derivative with
$M$ points and spatial step $\Delta x$ we approximate the derivative
$\partial_{xx} \approx A$  where
$A=\frac 1 {\Delta x^2} B \in \setR^{M-2, M-2}$
for Dirichlet boundary conditions and for Neumann boundary conditions,
 $$A=\frac 1 {\Delta x^2} \left(
   \begin{array}{rrrrr} 
 -2 & 2 \\     
    & B   \\
    & 2 & -2 
 \end{array}
 \right)\in \setR^{M, M}, \qquad \text{with} \qquad 
 B = \left( \begin{array}{ccccc}
     -2 & 1      & \\
     1 & \ddots & \ddots  \\
     & \ddots & \ddots & 1 \\
     &    & 1 & -2 \\
   \end{array}\right).$$
We choose not to use periodic boundary conditions since we compute 
a travelling front rather than a pulse and this would require 
a domain of twice the size and also introduces two travelling waves
that travel in opposite directions. 
Neumann and Dirichlet boundary conditions were shown to work well in
the deterministic case in \cite{thu05,thu06}. 
For the first spatial derivative 
we introduce
$(D_R u)_j = (u_{j+1} - u_j)/ {\Delta x}$, $(D_L u)_j = (u_{j} - u_{j-1})/ {\Delta x}$,
$(D_C u)_j = (u_{j+1} - u_{j-1})/ (2 \Delta x)$ for $j=1,...,M-1$ using 
Dirichlet boundary conditions 
$u^0 = \gamma_L, u_M=\gamma_R$ or Neumann boundary conditions
$u^0=u_1, u_{M-1}=u_M$. 
For convection terms we either use $D_C$ or, where
up-winding is an issue, we choose the appropriate $D_L$, $D_R$ or a
weighted combination \cite{bt04}
\begin{equation}
  \partial_x \approx D_h := \e{-\beta \mu} D_{L} + (1-\e{-\beta
    h}) D_{R},
\label{eq:Dbeta} 
\end{equation}
where $\beta$ is a parameter ($\beta = 0$ or $\beta = \frac 1 2$ in our
computations) and in what follows $h$ will be some function of the wave speed.
Recent work by Hairer and Voss \cite{HairerVoss} examine the discretization of
the advection term $u u_x$ for the stochastic Burger's equation and 
show that the limit is dependent on the discretization. The form of
advection for the Burger's equation and considered here is different
and as we can compare to cases without the advection term we did not
note and discretizatoin dependent 

Discretizing in space with $N$ grid points and after eliminating the
boundary conditions we obtain the following DAE system for $\lambda
\in \setR$ and $v\in\setR^{N-2}$ for  
Dirichlet or $v\in \setR^N$ for Neumann boundary conditions 
\begin{equation}\label{eq:dae}
  \begin{aligned}
     v' &= A v + \lambda (D_\lambda v + \eta) + f(v) + \varphi  \\
      0 & = \langle \uh_x, v-\uh \rangle,
  \end{aligned}
\end{equation}
where the vectors $\varphi, \eta$ are used to deal with the boundary conditions.
This system can be solved by using appropriate DAE solvers
\cite{ap98} or we can use a linear implicit Euler method to obtain the
fully discrete system 
\begin{equation}\label{eq:lin.euler}
  \begin{aligned}
    v^{n+1} &= v^n + \Delta t \left[ A v^{n+1} + \lambda^{n+1} (D_{\lambda^n} v^n + \eta) + f(v^n) + \varphi \right]\\ 
    0 & = \langle \Dc \uh, v^{n+1} -\uh \rangle 
  \end{aligned}
\end{equation}
which leads to
\begin{equation*}
  \begin{pmatrix}  
    I- \Delta t  A & - \Delta t  (D_{\lambda^n} v^n + \eta) \\
    \Delta x\, \Dc \uh^T & 0 
  \end{pmatrix}  
  \begin{pmatrix} v^{n+1} \\ \lambda^{n+1} \end{pmatrix}
  = \begin{pmatrix} v^n + \Delta t(f(v^n) + \varphi) \\ \langle \Dc \uh,\uh \rangle \end{pmatrix}.
\end{equation*}
Note that for the reference or template $\hat u_x$ we use the central
difference approximation $D_C$ since this is most accurate and 
convection instabilities are not an issue for this term. 

Under a uniqueness assumption of the travelling wave of the PDE 
it was shown in \cite{thu05} that for $L\to \infty$ and 
$\Delta x \to 0$ the stationary solution of \eqref{eq:dae} converges
to the exact travelling wave solution.
Moreover the solution of \eqref{eq:dae} inherits the
nonlinear stability properties of \eqref{eq:pde}. 
Numerically we observe below that the DAE system correctly computes the
travelling wave depending even when the travelling wave is not unique,
see \secref{sec:results}.

\subsection{Stochastic \PDE and stochastic travelling wave}
\label{sec:spdes}
We seek travelling wave solutions to 
the Stratonovich \SPDE
\begin{equation}\label{eq:spdeStrat}
  d u = \left[u_{xx} + f(u) \right]dt + g(u) \circ d W,\quad u(0)=u^0
\end{equation}
or 
the \Ito \SPDE
\begin{equation}\label{eq:spdeIto}
  d u = \left[u_{xx} + f(u) \right]dt + g(u) d W,\quad u(0)=u^0
\end{equation}
with  $g(u) = \nu + \mu h(u)$, where $\nu$
and $\mu$ are parameters that allow us to consider additive and
multiplicative noise. 
Results on the existence of a solution for
\eqref{eq:spdeStrat} and \eqref{eq:spdeIto} with $x\in\real$ 
domain can be found in \cite{Walsh86}. For the stochastic Allen-Cahn
equation with spatially smooth, bounded additive noise existence is
shown in \cite{Rgmnt02}. A recent paper \cite{Xie} shows existence for
the non-Lipschitz cases for space time white noise. 

We truncate the infinite domain and consider \eqref{eq:spdeStrat} (or \eqref{eq:spdeIto})
on a large finite domain so that $x\in[0,L]$ with either Neumann or
Dirichlet boundary conditions. 
For the finite domain with  $x\in [0,L]$  we refer to \cite[Theorem
7.4]{DaPrtoZbczyk} with $f$ and $g$ satisfying global Lipschitz
conditions, \cite{PrevotRoeckner} for some weaker conditions and
recent results in \cite{JentzenRoeckner} for non-Lipschitz $g$.
For the stochastic Nagumo equation we have $f(u)=u(1-u)(u-\alpha)$
and take $h(u)=u(1-u)$. 

We consider noise $W(t)$ to be a $Q$--Wiener
process~\cite{DaPrtoZbczyk}, and assume that covariance operator $Q$
and the linear operator $\partial_{xx}$ 
have the same eigenfunctions $\phi_j$. If the covariance operator has 
eigenvalues $\zeta_j\geq 0$ then we can write 
\begin{equation}
  \label{eq:W}
W(x,t) = \sum_{j\in\integers} \zeta_j^{1 / 2} \phi_j(x) \beta_j(t),  
\end{equation}
for independent Brownian motions $\beta_j$.
We take space-time noise that is white in time with exponential decay
in the spatial correlation length $\xi>0$ in which case  
$$\expect{dW(x,t)dW(y,s)} = C(x-y)\delta(t,s), 
\qquad C(x) = \frac 1 {2\xi} \exp\left(-\frac {\pi x^2}{4
    \xi^2}\right).$$
We approximate using \eqref{eq:W} by taking $\zeta_n= \exp(- \frac
{\xi^2 \lambda_j} {L})$, where $\lambda_j=\frac {j^2 \pi^2}{L^2}$  
and $L$ is the length of the interval \cite{Shardlow:05,Lrd+Rgmnt04}.

In the Stratonovich  case \eqref{eq:spdeStrat} we can eliminate 
the systematic effects of the noise on the drift and convert to an
\Ito integral.  This gives an
additional term to the nonlinearity so that the Stratonovich SPDE
\eqref{eq:spdeStrat} is equivalent to the \Ito SPDE 
\begin{equation}\label{eq:spdeCorrect}
  d u = \left[u_{xx} + \tilde{f}(u) \right]dt + g(u) d W
\end{equation}
where $\tilde{f}(u) = f(u) - C(0)g'(u)g(u)$, see for example
\cite{GrciaOjlvo+Sncho}.   
We can also convert from the \Ito interpretation \eqref{eq:spdeIto} to
a Stratonovich by 
\begin{equation}\label{eq:spdeCorrect2}
  d u = \left[u_{xx} + \tilde{f}(u) \right]dt + g(u)\circ d W
\end{equation}
where now $\tilde{f}(u) = f(u) + C(0)g'(u)g(u)$.

We discretize the SPDE in space by finite differences and evaluate
the noise on the spatial grid. 
In time we discretize with a constant time step $\Dt$. 
For the noise term we have an increment 
$$\Delta W_n= \sum_{j\in[-J,J]} \zeta_j^{1 / 2} \phi_j(x)
\xi_j,  $$ 
where $\xi_j \sim N(0,\Dt)$. 
To compute directly with the Stratonovich noise for
\eqref{eq:spdeStrat} we use the standard Heun method
\cite{Ojalvo,kloeden:1992} and also the semi-implicit Euler--Heun method 
\begin{equation}\label{eq:SPDEHeun}
  \begin{aligned}
    z & = u^n + g(u^n) \Delta W_n \\
    u^{n+1} &= u^n + \Delta t \left[ A u^{n+1}+ f(u^n) 
+\varphi \right] + \frac 1 2 (g(z) + g(u^n)) \Delta W_n
  \end{aligned}
\end{equation}
where $\Delta W_n$ in an increment of the noise and 
$\varphi$ arises from the boundary conditions. 

Although intuitively it is understood what is meant by a stochastic travelling
wave it is not easy to find a definition in the literature, however
for a review see \cite{GrciaOjlvo+Sncho,panja04}. Typically
a stochastic travelling wave and speed is either defined by the
evolution of a level set such as in \cite{Mllr_Swrs:95,Trbe96} or
through the evolution relative to a deterministic wave, such as
through a small noise expansion such as in \cite{Mkhlv_etal83}.
We will apply our methods to the case where in the deterministic case
the travelling wave is known to be unique and also where it is not unique.

Consider the \SPDE with a well defined wave with compact support as defined in
\cite{Shga94}, so that at 
$u(-\infty,\cdot\,)=u_-$, $u(\infty,\cdot\,)=u_+$. We can then define a
travelling wave and wave speed using the points
\begin{equation}
 a(t) := \sup\{z: u(x,t)=u_-, x\leq z\}, \quad
 b(t) := \sup\{z: u(x,t)=u_+, x\geq z\},
\label{eq:defab}
\end{equation}
and in  addition we can take the 'mid point' level set of a wave
\begin{equation}
 c(t) := \sup\{z: u(x,t)=(u_-+u_+)/2, x\leq z\}.
\label{eq:defc}
\end{equation}
These level sets define the position of the travelling wave. Note that
we take the supremum as there may be multiple crossings through the level
set (see for example \figref{fig:uhfail}) for the front. 
Given the positions an 'instantaneous' wave speed can be determined
from the $a$, $b$ and $c$  by differentiation.
We report below a wave speed $\Lambda_z(t)$, $z\in\{a,b,c\}$ defined by
\begin{equation}
\Lambda_z(t) = \E\left(\frac{z(t)-z(t_0)}{t-t_0}\right) \qquad z\in
\left\{a,b,c \right\} 
\label{eq:levLAM}
\end{equation}
where the expectation is taken over the number of realizations. We may
choose the initial time $t_0$ to either be at the start of the
computation $t_0=0$ or some later time ($t_0>0$) to avoid transient effects.
We believe This differs from the definition of wave speed used in
computations by \cite{Armro_etal,PhysRevE.58.5494,Mro:04} where they
report $z(t)/t$ for $z\in\left\{a,b,c\right\}$.
Numerically, these level set points $a,b,c$ are found by evolving 
the \SPDE \eqref{eq:spdeStrat} directly and interpolating over the grid.

If we assume the wave has some long time invariant speed an
alternative definition of the wave speed is to fit a linear polynomial
$P_{\LF}$ to the data  $(t,\E z(t))$, $z\in {a,b,c}$, $t\geq t_0$ where
\begin{equation}
P_{\LF}(t):=\LF t+K
\label{eq:LF}
\end{equation} 
Wave speeds may then be estimated by $\LF_z$, $z\in {a,b,c}$ where we may take
$t_0>0$ to avoid transients. Although this is a trivial extension of
wave speed defined by \eqref{eq:levLAM} we have not seen it reported
in the literature.

Finally we introduce a novel measure of the wave speed fpr SPDEs
through minimization of the $L^2$ norm $\|u(x,t)-\uh(x-y,t)\|^2$
against a fixed profile $\uh$. This is similar to the freezing the
wave in the deterministic case. We solve the SPDE and compute the
position $\gamma(t)$ of the wave. We 
then move the reference function relative to the travelling wave
solution $u$. That is we solve SPDE
$$  d u = \left[u_{xx} + f(u) \right]dt + g(u) \circ d W,\quad u(0)=u^0$$
and couple this to a reference function $\uh(x,t)$ that moves so that 
\begin{equation}
\label{eq:phase}
\langle \uh_x(t),u(t)-\uh(t) \rangle =0.
\end{equation}
Numerically this requires interpolation onto the spatial grid at each
time step. We compute an instantaneous wave speed $\lambda(t)$ and
this is related to the position of the wave through
$\gamma(t)=\int_0^t \lambda(s) ds$.
A wave speed $\LAM$ is then defined through the time average of
the instantaneous wave speed $\lambda$
\begin{equation}
\LAM(t) = \frac{1}{t-t_0}\int_{t_0}^{t} \E(\lambda(s)) ds .
\label{eq:LAM}
\end{equation}
So far we have not commented on the choice of profile $\uh$ to
minimize against. 

If a deterministic travelling wave profile exists then this is a
natural choice for $\uh$.
In examples where we do not have an analytic expression for the
deterministic travelling wave $\uh$ then this can be solved for
simultaneously or a sample solution solved for and saved.
However, this is a matter of choice and we could minimize the $L^2$ norm
against any fixed profile. 
One obvious choice of profile $\uh$ is to take the initial data, so
that $\uh=u(0)$.
Note that the choice of $\uh$ is important. In particular we
illustrate in \secref{sec:uref} that for a given $\uh$ the
minimization may not be unique and may fail numerically
if $\uh$ has small support.

We again note that from the position data $(t,\E \gamma(t))$ we can
fit a linear polynomial 
\begin{equation}
P_{\Lg}(t):=\Lg t +K
\label{eq:Lg}
\end{equation} 
to obtain an alternative estimate of the average wave speed.

\subsubsection{Freezing the stochastic travelling wave}
Inspired by the deterministic fixing of a wave we freeze a stochastic
travelling wave relative to a reference function $\uh$, this allows
direct computation of the wave speed from the $L^2$ minimization.
We take $\uh$ to be a fixed continuous function and we evaluate it
numerically at the grid pouints.
We examine the Stratonovich SPDE in a co-moving frame as we did for the
deterministic case.
First  let us examine the effects of a shift in space on the covariance of
the noise noting that
$$\E (dW(x+r(t),t),dW(y+r(s),s) ) = C(x+r(t)-y-r(s))\delta(t-s).$$
We see that for noise that is white in time the covariance is the
same in the two frames.

Let's consider the Stratonovich SPDE transforming to the moving frame
$u(x,t)=v(x+\gamma(t),t)$ we have
$$dv(x+\gamma(t),t) = v_x\circ d\gamma(t) + dv$$
and so formally we can write
$$
  \begin{aligned}
    d v &= \left[ v_{xx} + \frac{d\gamma(t)}{dt}v_x + f(v) \right] dt
    + g(v) \circ d W,\quad v(0)=u^0.
  \end{aligned}
$$
We introduce the phase condition to determine $d\gamma$ and 
minimize the $L^2$ norm, i.e. $\langle \uh_x, v-\uh \rangle=0$.
If we define the random variable $\lambda$ so that $\gamma(t) =
\int_0^t \lambda(s) ds$ then we seek to solve the SPDAE
\begin{equation}\label{eq:stratspdae}
  \begin{aligned}
    d v &= \left[ v_{xx} + \lambda v_x + f(v) \right] dt + g(v) \circ d W,\quad v(0)=u^0 \\
    0 & = \langle \uh_x, v-\uh \rangle
  \end{aligned}.
\end{equation}
In order to compute directly with the Stratonovich noise for
\eqref{eq:spdeStrat} we use either the standard Heun method
\cite{Ojalvo,kloeden:1992} or the semi-implicit Euler--Heun method 
\begin{equation}\label{eq:Heun}
  \begin{aligned}
    z & = u^n + g(u^n) \Delta W_n \\
    u^{n+1} &= u^n + \Delta t \left[ A u^{n+1} + \lambda^{n} \left( D_{\lambda^n} u^n +\eta \right) + f(u^n) +\varphi \right] 
    + \frac 1 2 (g(z) + g(u^n)) \Delta W_n\\   0 & = \langle \uh_x, u^{n+1}-\uh \rangle. 
  \end{aligned}
\end{equation}
The algorithm yields an approximation $u^n$, $n=0,1,2,\ldots$ to
$u(n\Delta t)$ and $\lambda^n$, and approximation to $\lambda(s)$.
This gives us a numerical scheme for the SPDAE in which the stochastic
travelling wave is frozen.

%
%
We have a time-dependent
random variable $\lambda(t)$ that we call the instantaneous 
wave speed. Of more physical interest is the time average of this
quantity $\LAMF$ that we call the wave speed and report
\begin{equation}\label{eq:deflambda}
\LAMF = \frac{1}{t-t_0}\int_{t_0}^{t} \E \lambda(s) ds.
\end{equation}
The instantaneous wave speed $\lambda$ gives the position $\gamma(t)$
of the wave.
If we assume the wave has some long time invariant
speed this can be estimated from $\gamma(t)$ by 
fitting a linear polynomial 
\begin{equation}
P_{\LgF}(t):=\LgF t +K
\end{equation} 
to the data from freezing the wave $(t,\E \gamma(t))$, for $t>t_0\geq 0$. 

The literature on solving stochastic DAEs is in its infancy, however there
are some analytic and computational results mainly arising from
examining noise in circuit simulations, see for example
\cite{SchnDnk98,Wnklr01,Wnklr03,Wnklr04}. We are not aware of work on
the existence directly for  \SPDAE.

Computing a travelling wave through \eqref{eq:stratspdae} or
we introduce  the random variable $\lambda$ which is used to freeze the wave. 
We could however define a weaker
versions by taking statistics of $\lambda$. For example we can take
the time-averaged wave speed $\Lambda(t)$ for each realization
\begin{equation}\label{eq:spdae.Lam}
  \begin{aligned}
    d v &= \left[ v_{xx} + \LAM v_x + f(v) \right] dt +
    g(v) \circ d W,\quad v(0)=u^0 \\ 
    0 & = \psi(v, \lambda).
  \end{aligned}
\end{equation}
Other weaker forms of travelling wave solution are possible where 
the instantaneous wave speed $\lambda$ or time average wave speed
$\Lambda$ of an individual realization is replaced by its expectation over
realizations, for example
\begin{equation}\label{eq:spdae.E}
  \begin{aligned}
    d v &= \left[ v_{xx} + \expect{\lambda} v_x + f(v) \right] dt +
    g(v) \circ d W,\quad v(0)=u^0 \\ 
    0 & = \psi(v, \lambda);
  \end{aligned}
\end{equation}
and
\begin{equation}\label{eq:spdae.ELam}
  \begin{aligned}
    d v &= \left[ v_{xx} + \expect{\Lambda} v_x + f(v) \right] dt +
    g(v) \circ d W,\quad v(0)=u^0 \\ 
    0 & = \psi(v, \lambda).
  \end{aligned}
\end{equation}
Using the sample mean of $\lambda$ and $\Lambda$ for fixing 
we are essentially using a ``group velocity'' to fix the wave and as a
result the mean profile will contain a spread as each individual
realization is not fixed at the same point.
By taking these weaker notions of wave speed to freeze the wave 
we observe spreading of the front profiles, as discussed in
\cite{GrciaOjlvo+Sncho}.

\section{Results for the Nagumo Equation}\label{sec:results}

We compare the different estimates of the wave speed and apply the
technique of freezing the wave to the Nagumo equation 
\eqref{eq:snagumo} for both multiplicative and additive
space-time white noise. 
For the majority of our simulations we take $\Dx=0.1$, $\Dt=0.05$ and a
spatial domain of $L=500$ or $L=800$ with Neumann boundary conditions 
and integrate till $t=100$. We compute $100$ realizations simultaneously.
In our computations of wave speeds unless stated we take $t_0=t/2$ to
reduce transient effects and we drop the dependence of the computed
wave speeds on the time $t$ and so report
$\LAM$, $\LAMF$, $\Lg$, $\Lambda_z$ and $\LF_z$,
$z\in\{a,b,c\}$.

\subsection{Deterministic PDE}
\label{sec:resdet}
Before we examine the stochastic PDE we briefly examine deterministic
computations. We point out some
features of computing the travelling wave and speed by direct
simulation of the PDE \eqref{eq:pde} versus freezing and solving
the PDAE \eqref{eq:dae}. 
In particular we examine the regime where the travelling wave is not
unique and the theory of \cite{thu05} on freezing the deterministic case no
longer holds. 
For $\alpha=-0.25$ the asymptotic travelling wave and wave speed
depends on the leading profile data of the initial data $u^0$. We take
two initial profiles $u^0(x)=u_{k_0}(x)$ with $k_0=0.1 < k_*$ and
$k_0=1/\sqrt{2}>k_*$.
To compute the speed by minimization we present results with 
reference functions $\uh=u_{\kh}$ with $\kh=\sqrt{2}$.
Results with a reference function with $\kh=0.1$ are identical
(see also \secref{sec:uref} for comments on the choice of reference function). 

We show in \tabref{tab:detwspdsnone} 
wave speeds computed from direct simulation of the PDE. 
For the PDE $\LAM$ is computed by moving the profile $\uh$ at the computed
wave speed from the minimization using the condition \eqref{eq:phase}. 
There is good agreement between the computed wave speeds, although
$\LAM$ appears to have converged faster than the other measures of the
wave speed to theoretical value. 
When we freeze the wave in the computational domain and solve the PDAE
we see from \tabref{tab:detwspdsfix} that the wave is frozen 
(to single precision) since the level set positions of $a(t)$,$b(t)$
and $c(t)$ do not change and hence wave speeds $\LF_z$, $z\in {a,b,c}$
from fitting the linear polynomial are zero (to single precision).  
The wave speed $\LAM$  estimated from freezing the wave and the
minimisation the $L^2$ norm agrees with the wave speeds
computed from the PDE (and is in fact a better approximation to the
theoretical values). 


\begin{table}
\begin{center}
\begin{tabular}{l|c|c|ccc|ccc} 
$\kh=1/\sqrt{2}$ & Theory & $\LAM$ &$\Lambda_a$ & $\Lambda_b$ & $\Lambda_c$ & $\LF_a$
 & $\LF_b$ & $\LF_c$ \\ \hline
$k_0=1/\sqrt{2}$ & 1.06066 &
1.06047 &
1.06025 &  1.06027 &  1.06026 & 
1.06025 &  1.06026 &  1.06026 \\
$k_0=0.1$ & $\geq 2.6$ & 
2.59741 &
2.59690 &  2.59689 &  2.59689 & 
2.59689 &  2.59689 &  2.59689 
\end{tabular}
\end{center}
\caption{Different measures of the wave speed computed from solving
  the deterministic PDE. To compute $\LAM$ the profile $\uh$ travels
  to minimize the $L^2$ norm \eqref{eq:phase}. Estimates of wave speeds
  $\Lambda_z$, are from \eqref{eq:levLAM} the level set and $\LF_z$
  from \eqref{eq:LF} from fitting, $z\in {a,b,c}$ . We see these
  measures of the wave speed agree to $4$ decimal places.
}
\label{tab:detwspdsnone}
\end{table}

\begin{table}
\begin{center}
\begin{tabular}{l|c|c|ccc|ccc} 
$\kh=1/\sqrt{2}$ &  Theory & $\LAMF$ & $\Lambda_a$ & $\Lambda_b$ & $\Lambda_c$ & $\LF_a$
 & $\LF_b$ & $\LF_c$ \\ \hline
$k_0=1/\sqrt{2}$ & 1.06066 &
1.06052  &
2.0e-07 &-1.1e-06 &1.5e-08&
1.8e-07 &-9.5e-07 &1.4e-08\\
$k_0=0.1$ & $\geq 2.6$ & 
2.60048  &
3.1e-06 &-3.6e-07 &-1.4e-07&
1.3e-06 &-8.4e-08 &-6.0e-08 
\end{tabular}
\end{center}
\caption{Wave speeds computed from solving the PDAE and freezing the
  travelling wave. The fact that the wave does not move in the domain
  can be seen from the level set wave speeds $\Lambda_z$ and $\LF_z$,
  $z\in\{a,b,c\}$  which are close to zero. The wave speed  $\LAMF$
  agrees with that computed for the PDE given in \tabref{tab:detwspdsnone}.}
\label{tab:detwspdsfix}
\end{table}

\subsection{Stochastic travelling wave and frozen wave}
\begin{figure}[hth]
\begin{center}
  (a) \hspace{0.32\textwidth} (b)  \hspace{0.32\textwidth} (c)  \\
  \includegraphics*[width=0.32\textwidth]{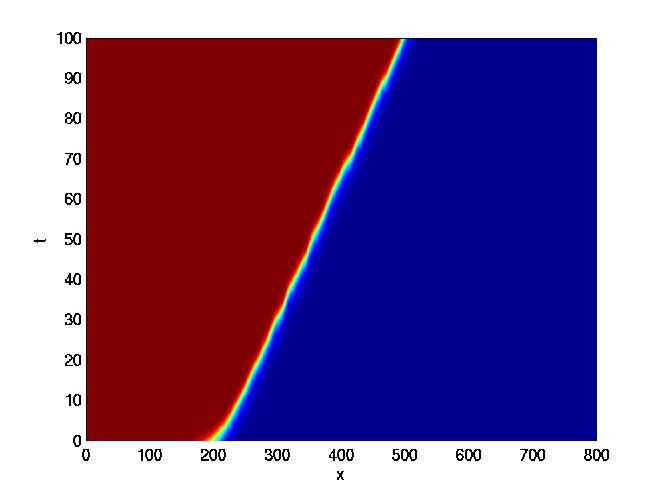}
  \includegraphics*[width=0.32\textwidth]{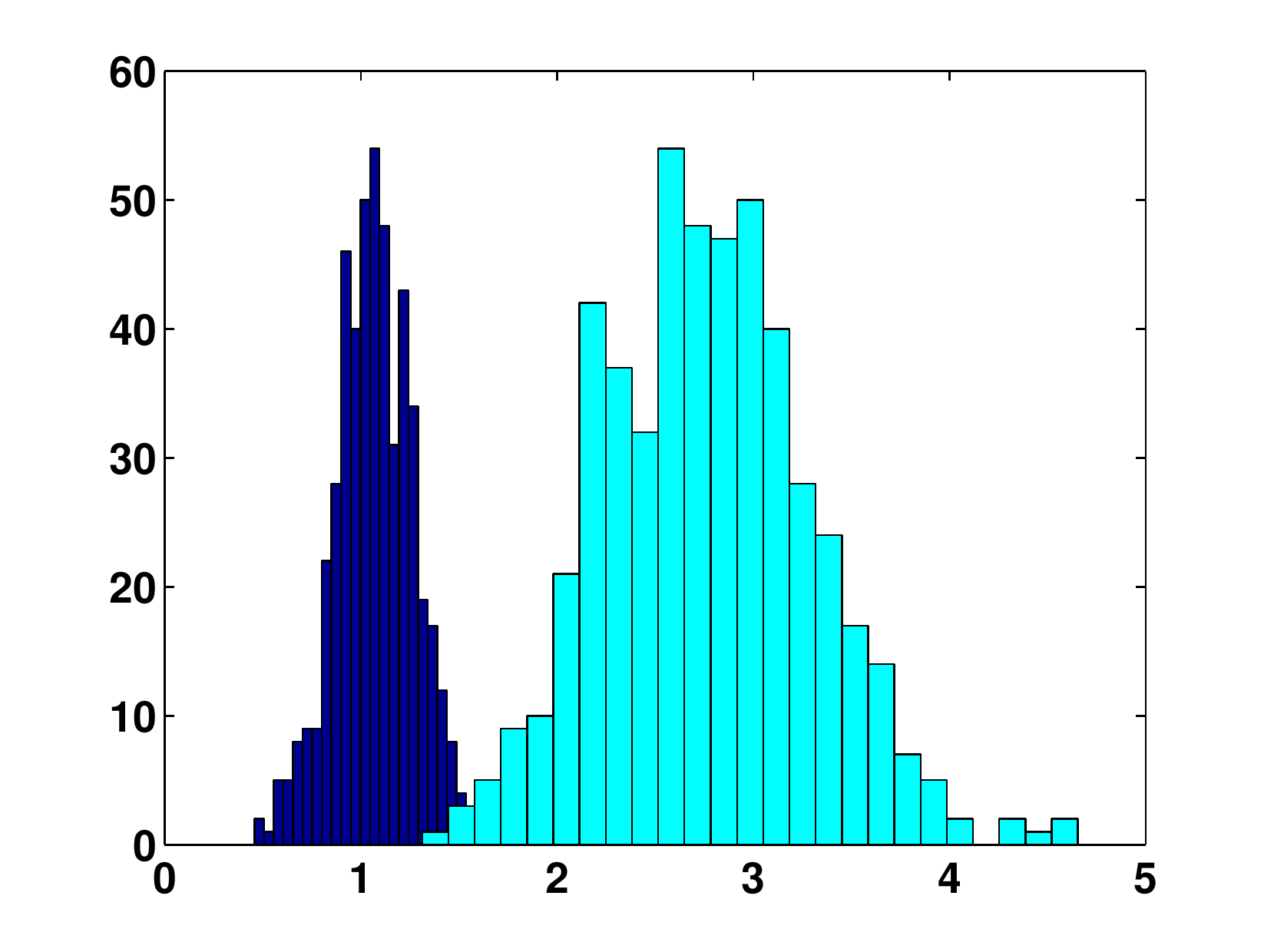}
  \includegraphics*[width=0.32\textwidth]{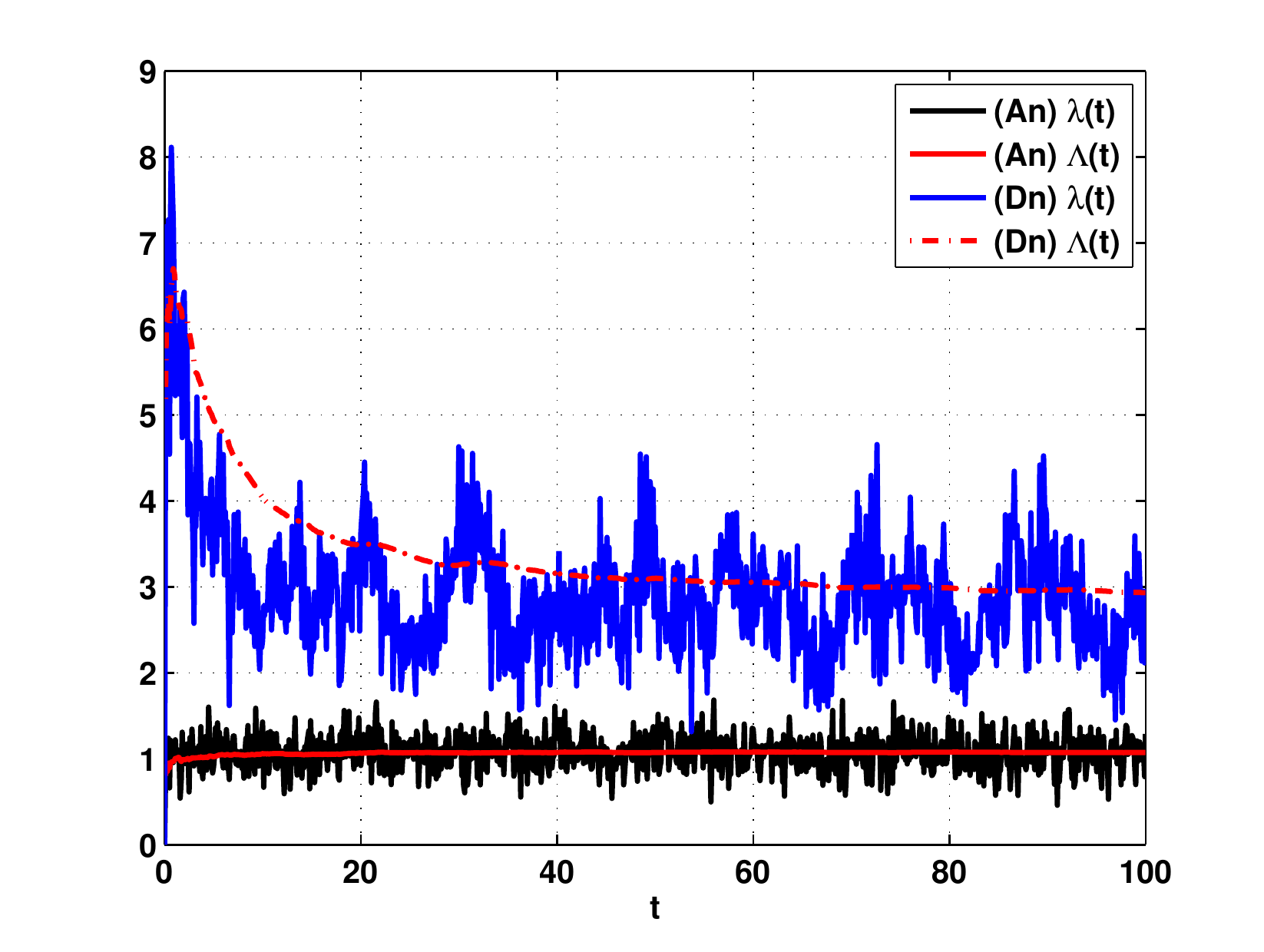}\\
  \caption{(a) Space-time plot of a single realization of the SPDE
    showing a travelling wave. In (b) distributions of the
    instantaneous wave speeds $\lambda$ computed for the two different
  initial data sets. In (c) we plot $\lambda(t)$ and the time averages
  $\LAM(t)$ with $t_0=0$, $t_1=t$.}
  \label{fig:nag_none_a-025}
\end{center}
\end{figure}

To illustrate computations for the stochastic PDE 
we start by taking Stratonovich multiplicative noise with $\mu=0.1$ and a
correlation length of $\xi=0.1$.  In \tabref{tab:wspdsnone1} we show
results from solving the SPDE with the same single realization of the
noise using the same two 
different sets of initial data and two difference reference functions
as for the deterministic case. 
Since we have taken the same noise realization, when initial data is the
same our measures of the wave speed $\Lambda_z$, and $\LF_z$, $z\in
{a,b,c}$ are identical and independent of the reference function $\uh$. 
The choice of reference function does change
the wave speed measured by the minimization in approximately the fourth
decimal place (compare $\LAM$ or $\Lg$ for the different $\kh$
values in \tabref{tab:wspdsnone1}). This small difference is due to a
combination of interpolation errors and is not seen for the
SPDAE below where we do not need this interpolation.
If we compare the values of the deterministic PDE \tabref{tab:detwspdsnone}
and SPDE case \tabref{tab:wspdsnone1} for the single realization we see the
wave speeds with noise are slightly larger than the deterministic
case. 
In \figref{fig:nag_none_a-025} we plot the result of a
single realization in (a) for the SPDE with initial data $u^0$ with
$k_0=1/\sqrt{2}$. The wave front is initially at $x\approx 200$ and 
travels to $x\approx 500$. For the two different initial data we have
plotted in (b) the two distributions of the instantaneous wave speed 
$\lambda$ used to compute the wave speed through the
minimization (with $\kh=1/\sqrt{2}$). The mean of these distributions
gives the corresponding wave speeds, $1.07575$ for
$k_0=1/\sqrt{2}$ and $2.77146$ for $k_0=0.1$. For initial data
$k_0=1/\sqrt{2}$ with a wave speed of $1.07575$ the variance of
$\lambda$ is smaller.
In (c) we have plotted for the two different initial data sets 
the instantaneous wave speed $\lambda(t)$ and the corresponding time
averaged wave speeds $\LAM(t)$ with $t_0=0$, $t_1=t$. We see faster
convergence of the wave speed for initial data $k_0=1/\sqrt{2}$ and
again the reduced variability in the instantaneous wave speed.

\begin{table}
\begin{center}
\begin{tabular}{l|cc|ccc|ccc} 
$t_0=50, t_1=100$ & $\LAM$ & $\Lg$ & $\Lambda_a$ & $\Lambda_b$ &
$\Lambda_c$ & $\LF_a$ & $\LF_b$ & $\LF_c$ \\ 
\hline
{\small $k_0=1/\sqrt{2}, \kh=1/\sqrt{2}$}& 
1.07575 &  1.06965 &
1.07484 &  1.06746 &  1.07536 & 
1.07149 &  1.07147 &  1.06935 \\
{\small $k_0=1/\sqrt{2}, \kh=0.1$} & 
1.07538 &  1.06987 &
1.07484 &  1.06746 &  1.07536 & 
1.07149 &  1.07147 &  1.06935 \\
{\small $k_0=0.1, \kh=0.1$} & 
2.77250 &  2.78015 &
2.80594 &  2.74516 &  2.76020 & 
2.79657 &  2.72070 &  2.78112 \\
{\small $k_0=0.1, \kh=1/\sqrt{2}$} & 
2.77146 &  2.78178 &
2.80594 &  2.74516 &  2.76020 & 
2.79657 &  2.72070 &  2.78112 \\
\end{tabular}
\end{center}
\caption{Wave speeds computed from solving a single realization of the
  SPDE with noise  intensity $\mu=0.1$ and correlation length $\xi=0.1$. To
  compute $\LAM$ the profile $\uh$ travels with the appropriate speed
  found by minimization of the $L^2$ norm.}
\label{tab:wspdsnone1}
\end{table}

Let us now compare to a single realization where the stochastic
travelling wave is frozen and we solve the SPDAE \eqref{eq:stratspdae}
using the Heun method \eqref{eq:Heun}. 
Results on the wave speeds for the two initial data sets and two
reference functions are reported in \tabref{tab:wspdsfix1}. 
We chose $\Lambda_c$ and $\LF_c$ to represent values computed form the
level set approaches. The level sets no longer travel (on average) and
hence have wave speeds with values close to zero. 
Note that the noise path is not the same as solving the SPDE for 
\tabref{tab:wspdsnone1} and so we do not expect the values to be
exactly the same, they are however close.
The wave speed $\LAMF$ estimated by the minimization is identical for
the two different reference functions (solving the SPDAE we do not
have the same interpolation errors as when solving the SPDE). 
However, the choice of reference function is an issue for the SPDE and
we consider this further in \secref{sec:uref}.

\begin{figure}[hth]
\begin{center}
     (a) \hspace{0.32\textwidth} (b)  \hspace{0.32\textwidth} (c)  \\
  \includegraphics*[width=0.32\textwidth]{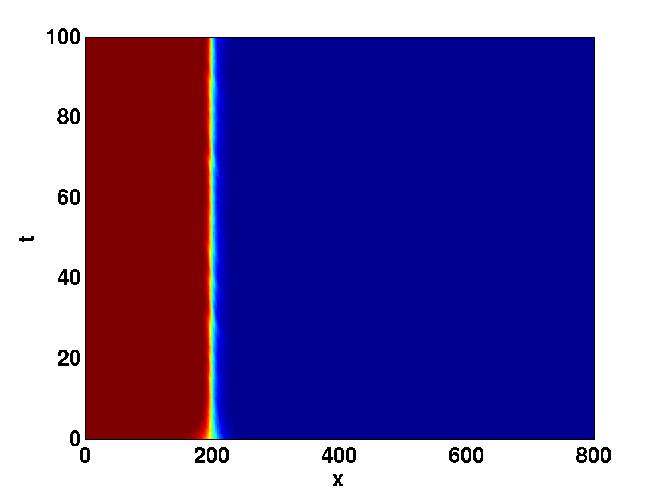}
  \includegraphics*[width=0.32\textwidth]{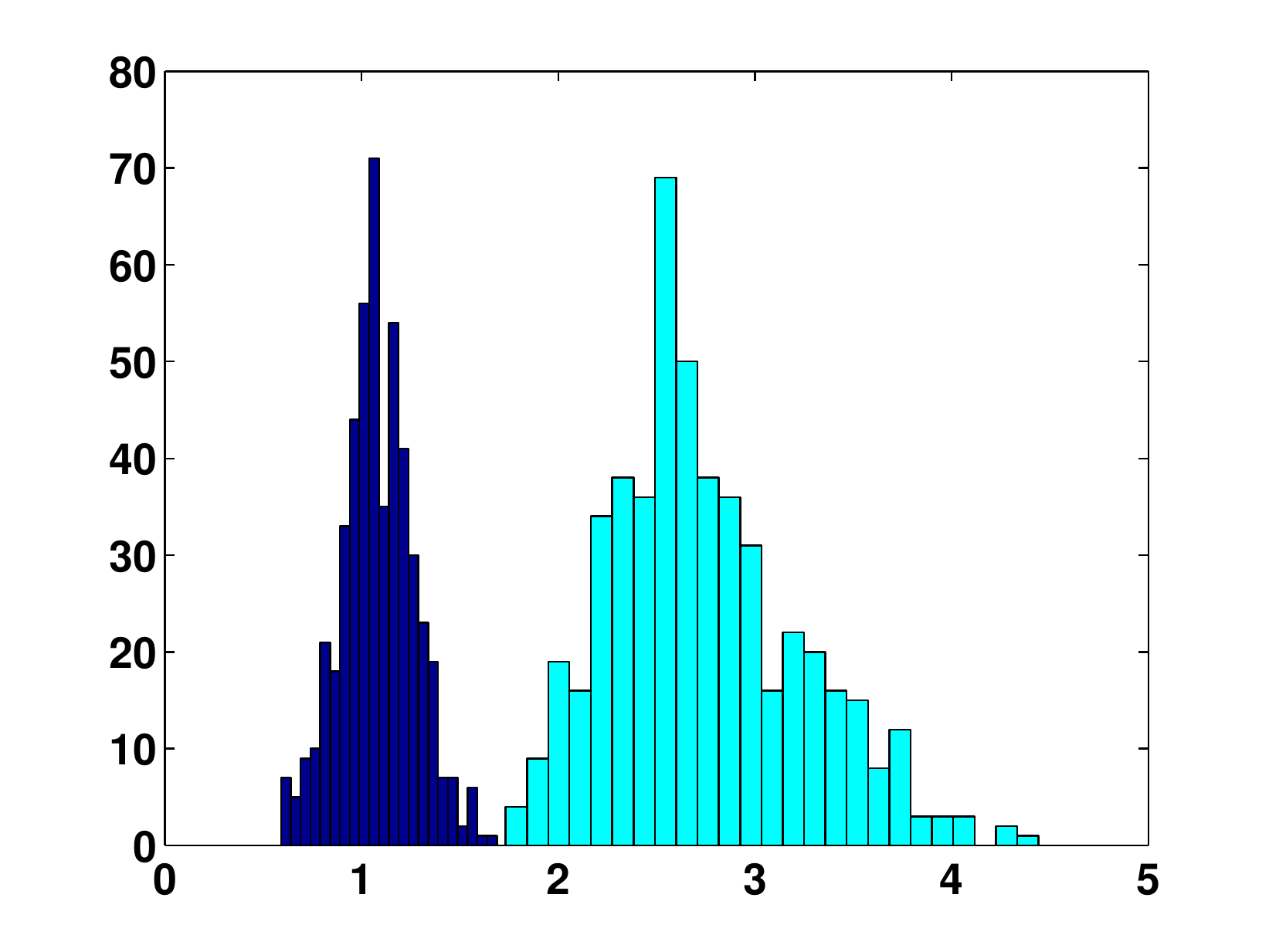}
  \includegraphics*[width=0.32\textwidth]{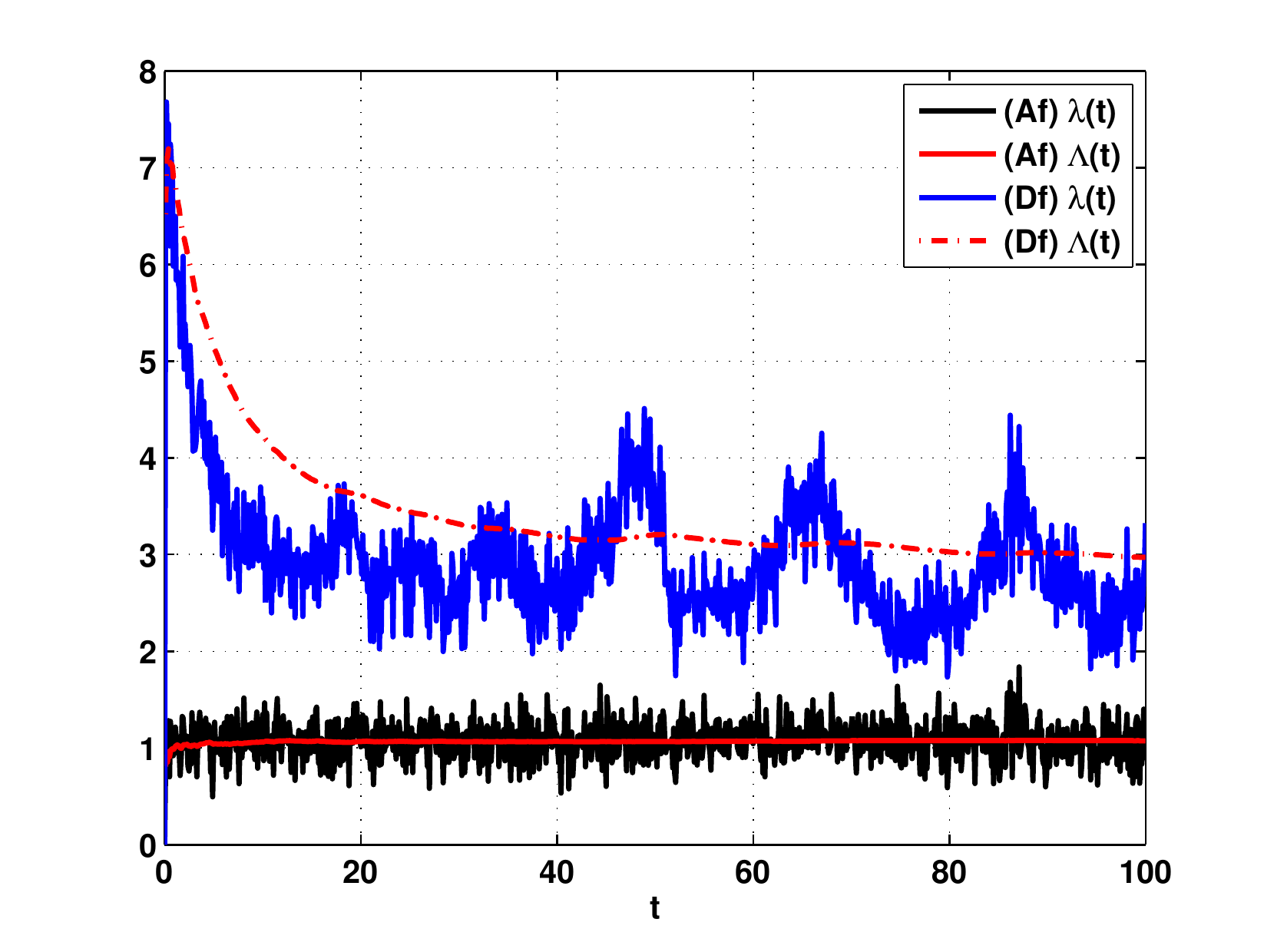}
  \caption{(a) Space-time plot of a single realization of the frozen SPDE
    showing a travelling wave. In (b) distributions of the
    instantaneous wave speeds $\lambda$ computed for the two different
  initial data sets. In (c) we plot $\lambda(t)$ and the time averages
  $\LAM(t)$ with $t_0=0$, $t_1=t$. }
  \label{fig:nag_fix_a-025}
\end{center}
\end{figure}

In \figref{fig:nag_fix_a-025} we have plotted in (a) the space-time
plot of solution of the SPDAE. The front starts at $x\approx 200$ and
remains (on average) at that position throughout the computation
illustrating that the wave does not travel (compare to
\figref{fig:nag_none_a-025} (a)). In   
(b) for the two different initial data $k_0=1/\sqrt{2}$ (with mean
$1.08522$) and $k_0=0.1$ (with mean $2.74311$) we
have plotted the two distributions of  the instantaneous wave speed
$\lambda$ used to compute the wave speed through the minimization
(with $\kh=1/\sqrt{2}$). Comparing with  
\figref{fig:nag_none_a-025} (b) we see similar distributions and 
greater variance with initial data with $k_0=0.1$ than
$k_0=1/\sqrt{2}$ (as in \figref{fig:nag_none_a-025}).
In (c) we have plotted for the two different initial data sets 
the instantaneous wave speed $\lambda(t)$ and the time averaged wave
speed $\LAM(t)$ with $t_0=0$. We see faster convergence of the wave
speed $\Lambda(t)$ for $k_0=1/\sqrt{2}$ than for  $k_0=0.1$.
%
%
\begin{table}
\begin{center}
\begin{tabular}{l|cc|ccc|ccc|cc} 
$t_0=50$, $t_1=100$ & $\LAMF$ & $\LgF$  & $\Lambda_c$  & $\LF_c$ \\ \hline
$k_0=1/\sqrt{2}, \kh=1/\sqrt{2}$ & 
1.08522   & 1.08828 &
2.421e-04& -2.131e-04\\

$k_0=1/\sqrt{2}, \kh=0.1$ & 
1.08522 & 1.08828 &
2.421e-04& -2.131e-04 \\
$k_0=0.1, \kh=0.1$ & 
2.74311 & 2.76005 &
-5.703e-02& -6.404e-03\\
$k_0=0.1, \kh=1/\sqrt{2}$ & 
2.74311 & 2.76005 & 
-5.703e-02& -6.404e-03\\
\end{tabular}
\end{center}
\caption{Wave speeds computed from solving a single realization of the
  SPDAE with noise intensity $\mu=0.1$ and correlation length $\xi=0.1$.} 
\label{tab:wspdsfix1}
\end{table}

Rather than looking at a single realization more physically meaningful
results are found from taking the expectation over many realizations. In
\tabref{tab:wspds10100} we examine wave speeds based on $100$
realizations of both the SPDE and SPDAE. The different measures of the
wave speed are in broad agreement. The larger
uncertainty in $\LAM$ and $\LAMF$ originates in the large variance in the
instantaneous wave speeds $\lambda$ and is a drawback of the
minimization approach.

We also compare the profiles from the \SPDE to profiles obtained 
from the \SPDAE. To avoid the spreading of the wave we need to align
individual realizations of the \SPDE. We chose as a common reference
the level set $c(100)$. If we examine the final time profiles for the
runs we find that the weak error, 
$\|\expect{u_{SPDAE}(100)} - \expect{u_{SPDE}(100)}\|_{L^2}^2$ 
for $10$ realizations is $\approx 0.0150$ and for $100$ realizations
$\approx 0.0144$ and $\approx 0.0117$ with $1000$ realizations.

\begin{table}
\begin{center}
\begin{tabular}{l|c|c|cc} 
&  $\LAM$ or $\LAMF$ & $\Lg$ or $\LgF$ & $\Lambda_c$ &  $\LF_c$ \\ \hline
SPDE   & 
1.08588 $\pm$ 0.19680  & 1.08388 $\pm$ 2.73e-03  &
1.08381 $\pm$2.81e-03  & 1.08390 $\pm$2.68e-03 \\ 
SPDAE &
1.08951 $\pm$ 0.19512 &  1.08790 $\pm$ 2.39e-03 & 
-4.0e-05 $\pm$ 2.0e-5 & 3.0e-05  $\pm$ 2.0e-5
\end{tabular}
\end{center}
\caption{Expected values of the wave speeds taken over 100
  realizations solving the SPDE and the SPDAE. Initial data taken with
  $k_0=1/\sqrt{2}$,and reference function with $\kh=1/\sqrt{2}$.}
\label{tab:wspds10100}
\end{table}

\begin{figure}[hbt]
\begin{center}
  \psfrag{LAM}{$\LAM$}
  \psfrag{mu^2}{$\mu^2$}
     (a) \hspace{0.45\textwidth} (b)  \\
  \includegraphics*[width=0.42\textwidth]{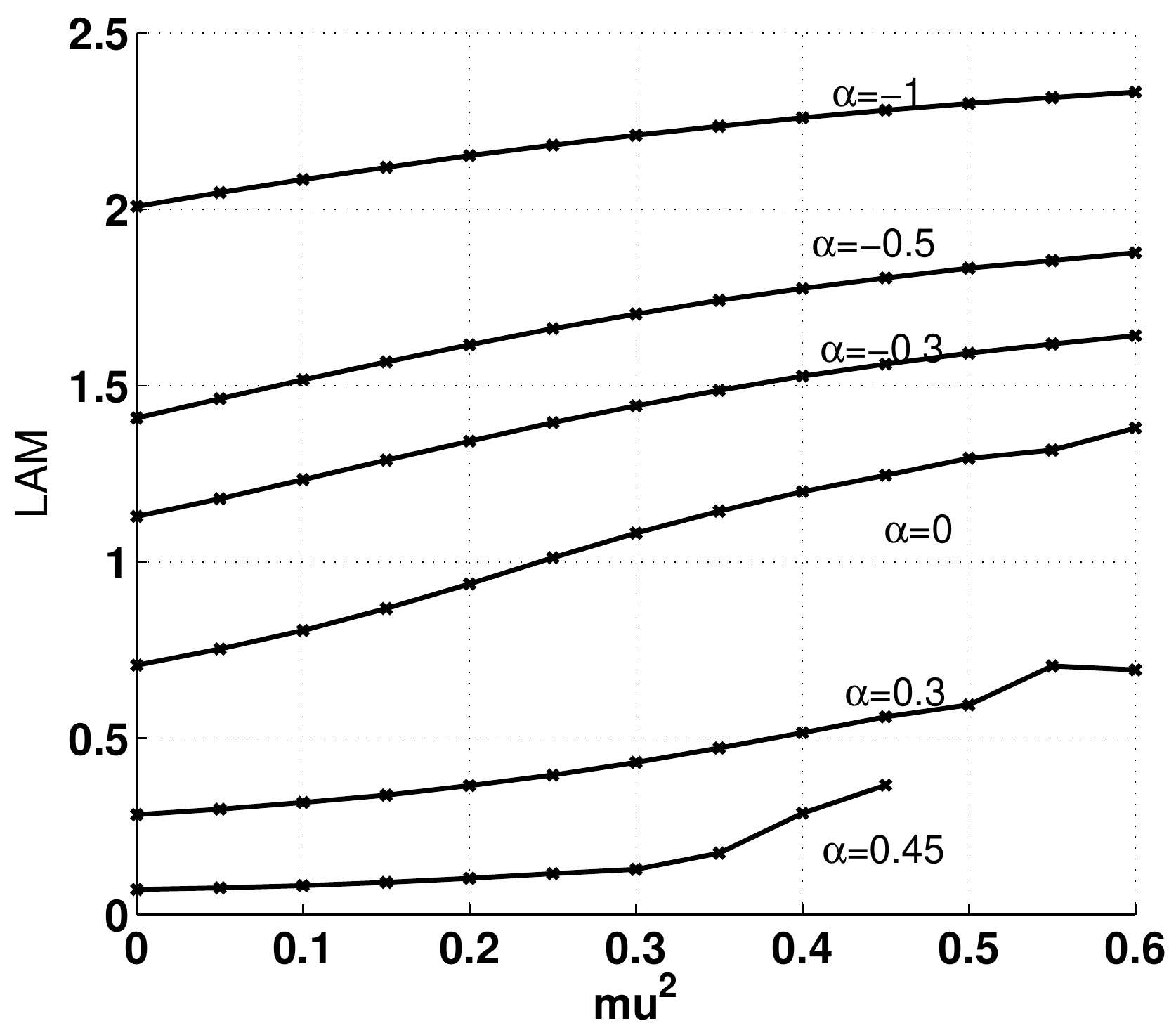}
  \includegraphics*[width=0.45\textwidth]{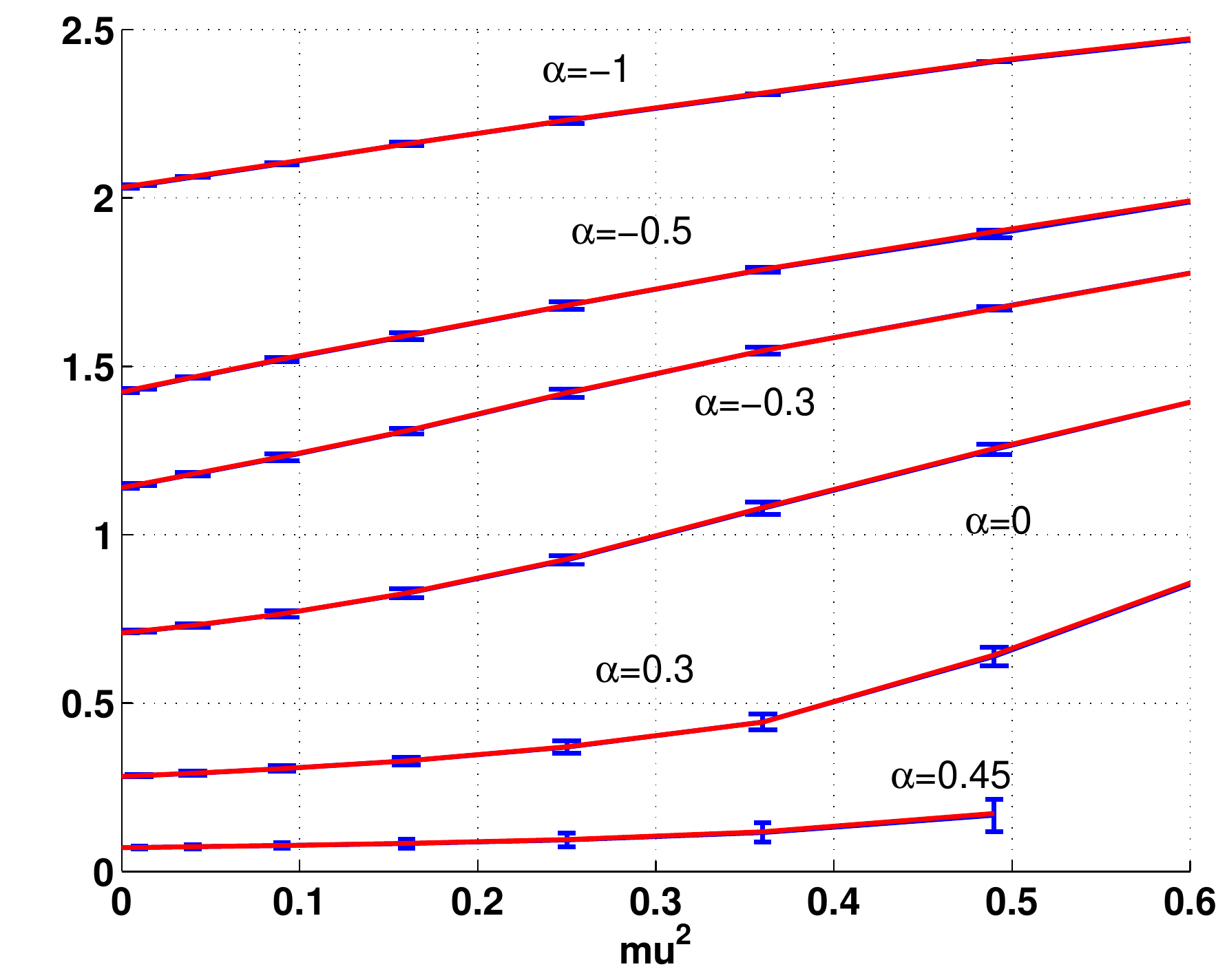}\\
  \caption{Wave speeds $\LAM$ with increasing noise intensity for
    Stratonovich noise with correlation length equal to that of the
    grid. In (a) solving the SPDAE where the wave is frozen with
    $\xi=0.5=\Dx$ (b) the SPDE with $\xi=0.1=\Dx$. Each line
    corresponds 
    to a different nonlinearity with
    $\alpha\in\{-1,-0.5,-0.3,0,0.3,0.45\}$.}
\label{fig:SPDE_SPDAE_compare}
\end{center}
\end{figure}

In \figref{fig:SPDE_SPDAE_compare} we compare results for the SPDAE (a) and
SPDE (b) for a range of different nonlinearity's $\alpha\in\{-1,-0.5,-0.3,0,0.3,0.45\}$ and noise intensities 
measured by $\mu^2 \in [0,1]$. Initial data approximates a
step function and the spatial correlation length of the noise $\xi$ is that
of the computational grid $\Dx$.  The results in (a), where
$\Dx=\xi=0.5$ where agree with those in 
\cite{Armro_etal,PhysRevE.58.5494}, reproduced in  
\cite{GrciaOjlvo+Sncho}, where the authors obtain a 
front velocity taking an average over an ``appropriate time window''
of $\int_L u(x,t) dt$ and compare to a small noise analysis. In (b) we
took a smaller spatial step $\Dx=0.1$ and have plotted the wave speed $\LAM$
computed both from minimization and from the level set,
$\Lambda_c$ on which the error bars are based.

We observe that
the effect of the two different approximations to spatially white
noise is to increase the speed of the wave.
Note that some realizations where the wave is frozen in
\figref{fig:SPDE_SPDAE_compare} (a) fail to exist due to numerical
instability, see \secref{sec:spdae_unstable}.
%

%

\subsubsection{Numerical instability}
\label{sec:spdae_unstable}
The numerical approximation of SDEs and SPDES where the solution is  
constrained in phase space is an area under development. 
For the Nagumo equations \eqref{eq:spdeStrat} (or \eqref{eq:spdeIto})
$u\in[0,1]$, another typical example is a positivity constraint where
$u>0$. Numerical instability can lead to non-physical solutions and
potentially to unphysical unbounded growth of the numerical solution.
A number of approaches have been proposed to 
simulations to enforce constraints on the numerics and a review of
these types of methods for SDEs is contained in \cite{LordKoekkoekvanDijk}.  
One method to avoid unbounded growth in numerics from nonphysical
solutions the nonlinearity and noise can be adapted as in 
\cite{MoroSchurz07,Drng+etal:03,Drng.etal:05,Shardlow:05}.

\begin{figure}[hth]
\begin{center}
  (a) \hspace{0.48\textwidth} (b)  \\
  \includegraphics*[width=0.48\textwidth]{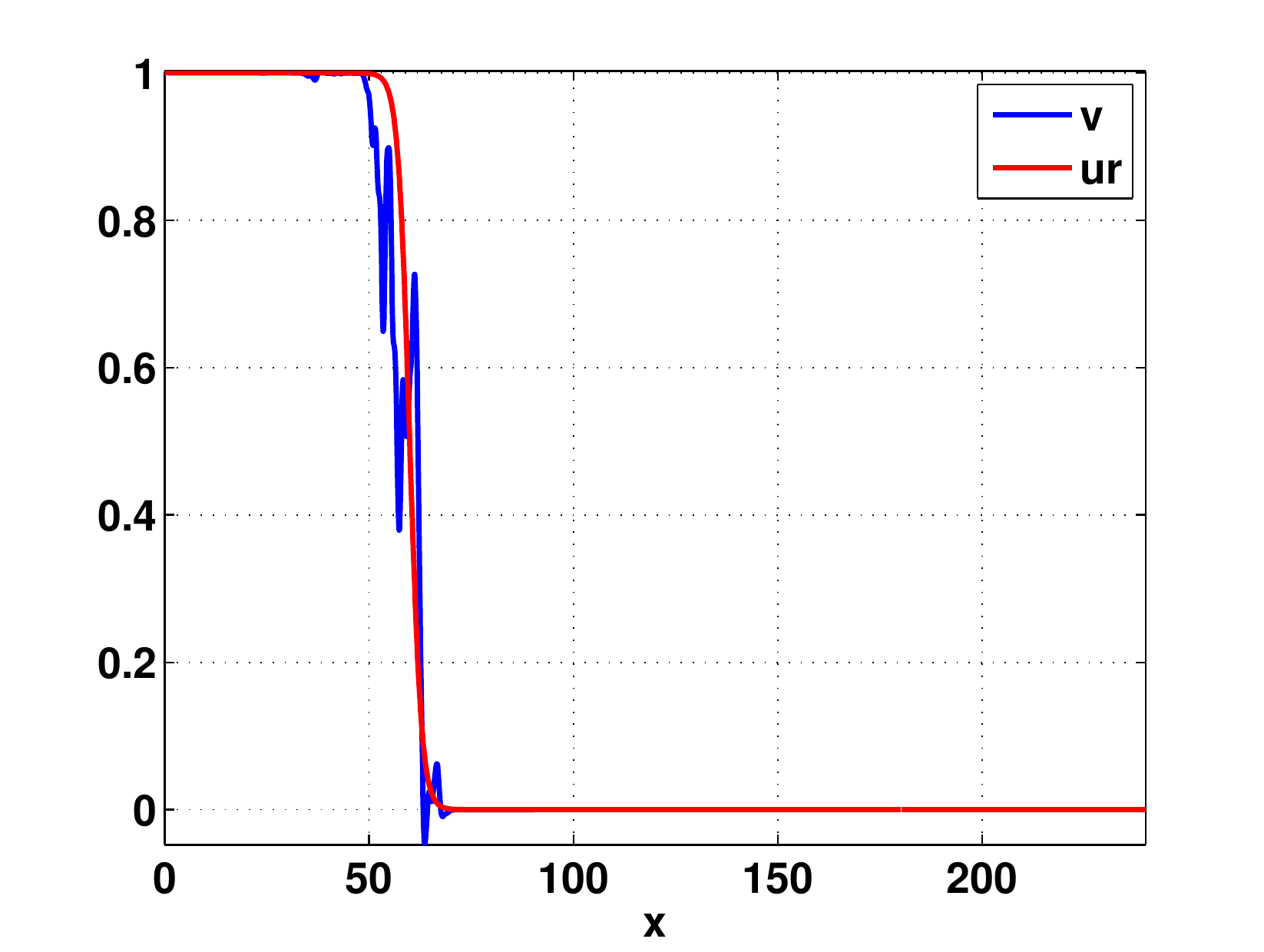}
  \includegraphics*[width=0.48\textwidth]{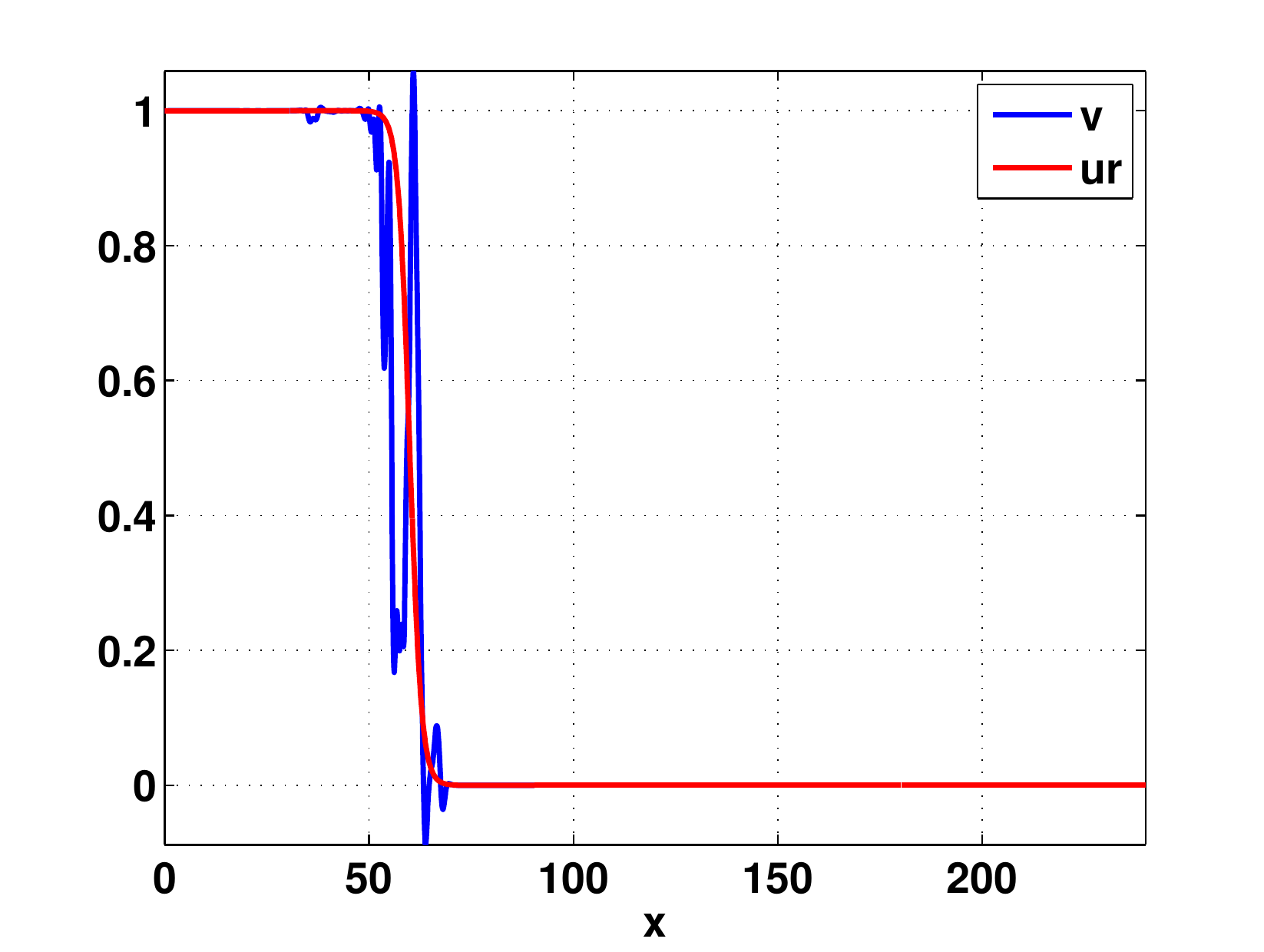}\\
  \caption{Plot of realization of a single realization of the noise
    illustrating instability (a) at $t=34.7$ and (b) $t=35.2$, with
    corresponding instantaneous wave speeds of $10.98$ and $-10.60$.}
  \label{fig:instab}
\end{center}
\end{figure}

We found that solving the SPDEs \eqref{eq:spdeStrat} (or
\eqref{eq:spdeIto}) such instability was not an issue. However, when
freezing the wave and solving the SPDAE \eqref{eq:stratspdae} did lead 
to non-physical solutions. In \figref{fig:instab} we have
frozen the wave and show one realization at $t=34.7$ (a) with an
instantaneous wave speed $10.98$ and (b) $t=35.2$ with an
instantaneous speed $-10.60$. The non-physical regions where $u<0$ and
$u>1$ then grow in magnitude with further iterations.

\figref{fig:nag_none_a-025} and \figref{fig:nag_fix_a-025} show the
distribution of the instantaneous wave speed $\lambda$ as computed for
the SPDE and SPDAE respectively where the wave is frozen. 
The SPDAE system \eqref{eq:stratspdae} includes the advection term
$\lambda v_x$, where $\lambda$ is a random variable with a particular 
distribution - which may lead to either large positive and/or negative 
values of $\lambda$. 
Numerically this is particularly true for large
noise intensities or small correlation lengths of the noise.
The result of this is a loss of numerical stability.
Although we were able to control unbounded growth by modifying the
equation solved close to the $u=0$ and $u=1$, direct comparisons
of wave speeds to SPDE calculations showed this can lead to a bias
in the estimate, so we do not include such results here.
Hence, in \figref{fig:SPDE_SPDAE_compare} (a) results for the SPDAE
equation are reported with the expectation taken over solutions that
existed to the final time. A large number of initial realizations was
taken so that the final expectation is over at least $1000$
realizations. Although we observe the same results calculating the
wave speed based on level set methods, in general not taking the results where
there is numerical blow up may bias the statistics.

%
%
%

\subsubsection{Choice of reference function $\uh$ for the minimization}
\label{sec:uref}

We commented early that natural choice for the reference function
would be to take either a deterministic travelling wave or the initial
data. However, the width of the reference function $\uh$ plays an
important role in the computed wave speed for both the SPDE and
SPDAE. If we take the reference function $\uh$ to be the Heaviside
function then the minimization of the $L^2$ norm fails.
We observe numerically that narrow reference functions can also lead
to numerical failure of the minimization and indeed there may be more
than one minimal position. 

\begin{figure}[hth]
\psfrag{lm}[][][0.6]{$\lambda$}
  \begin{center}
    (a) \hspace{0.48\textwidth} (b)  \\
    \includegraphics*[width=0.42\textwidth]{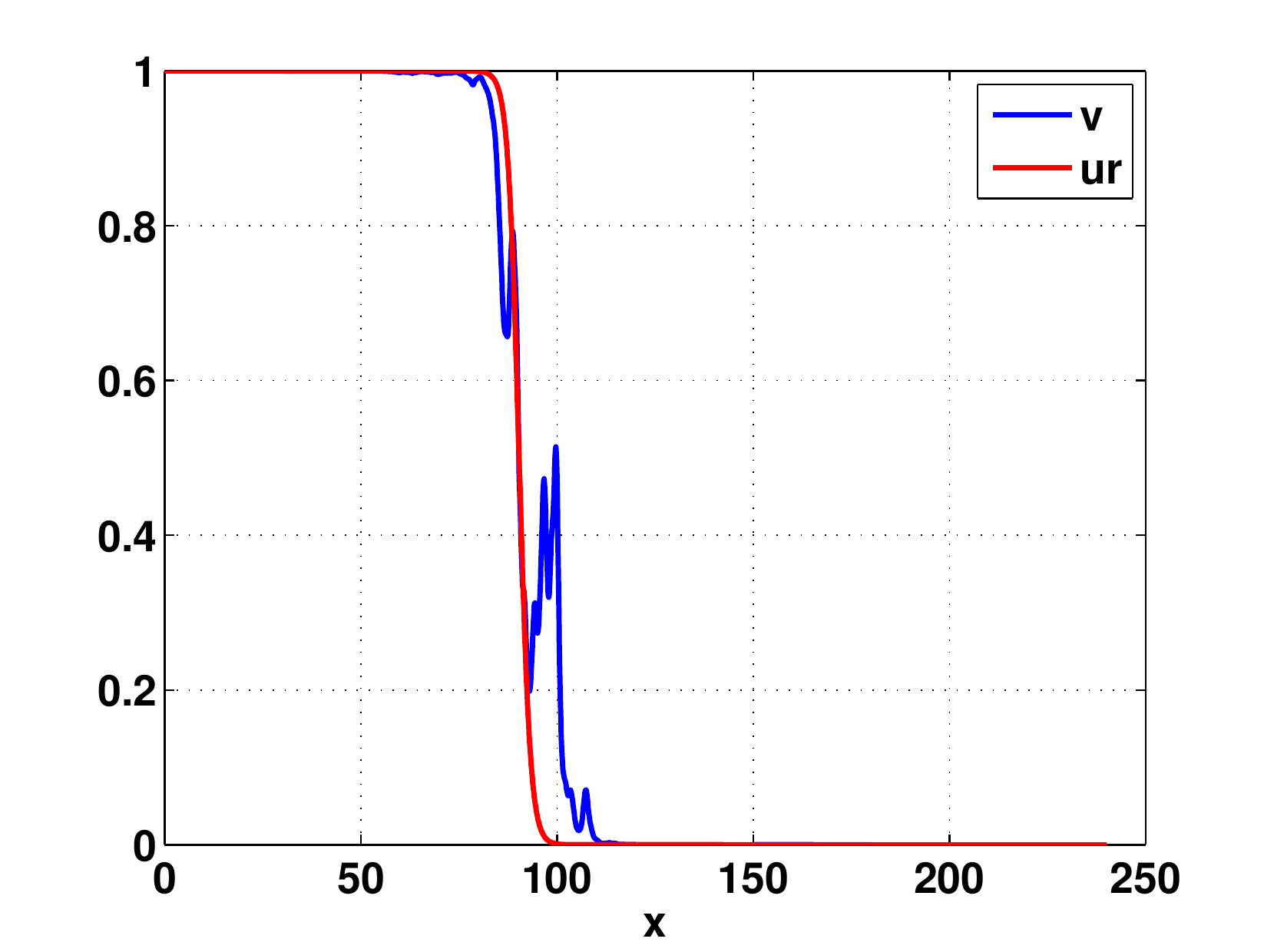} 
    \includegraphics*[width=0.4\textwidth]{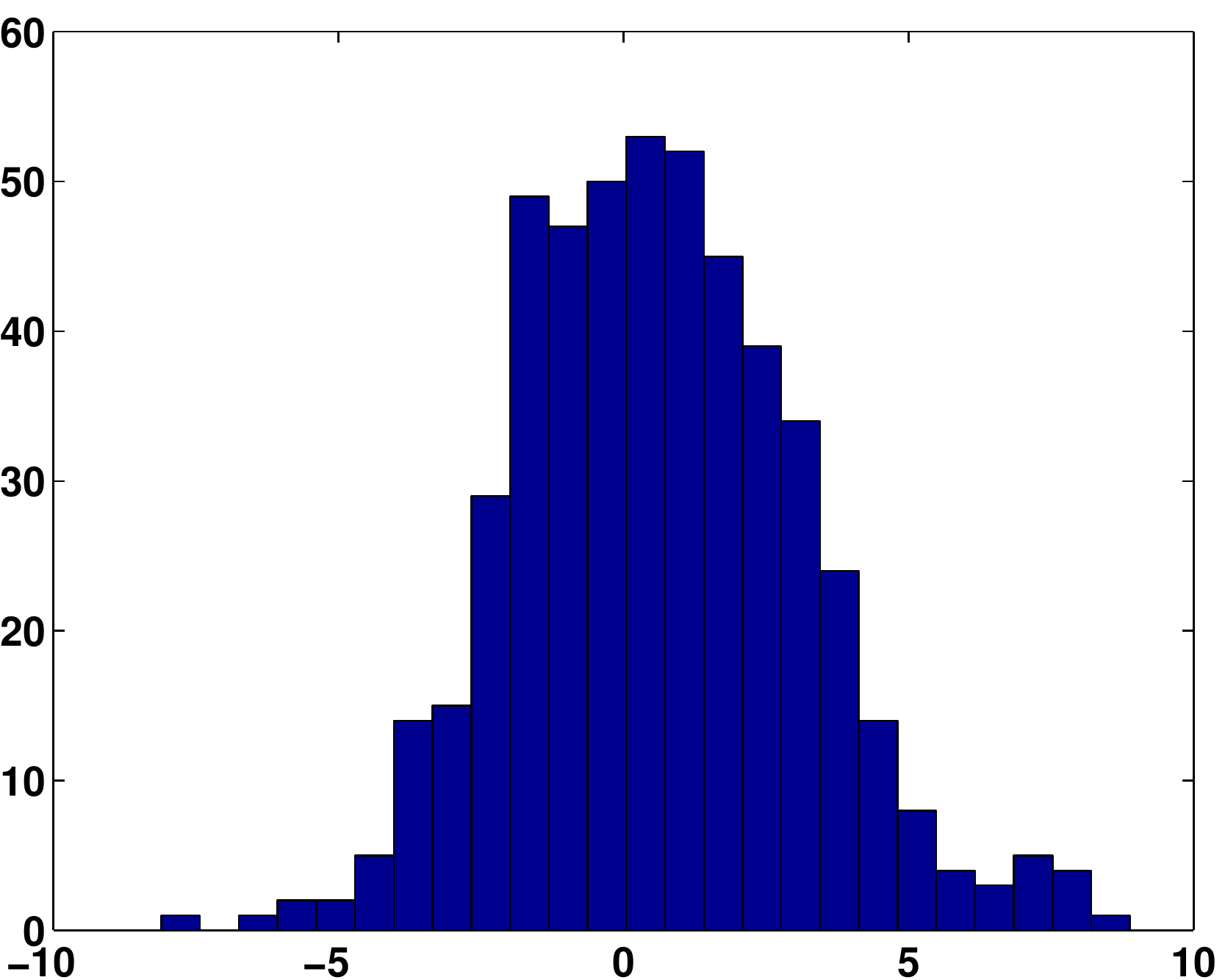} \\
    (c) \hspace{0.48\textwidth} (d)  \\
    \includegraphics*[width=0.42\textwidth]{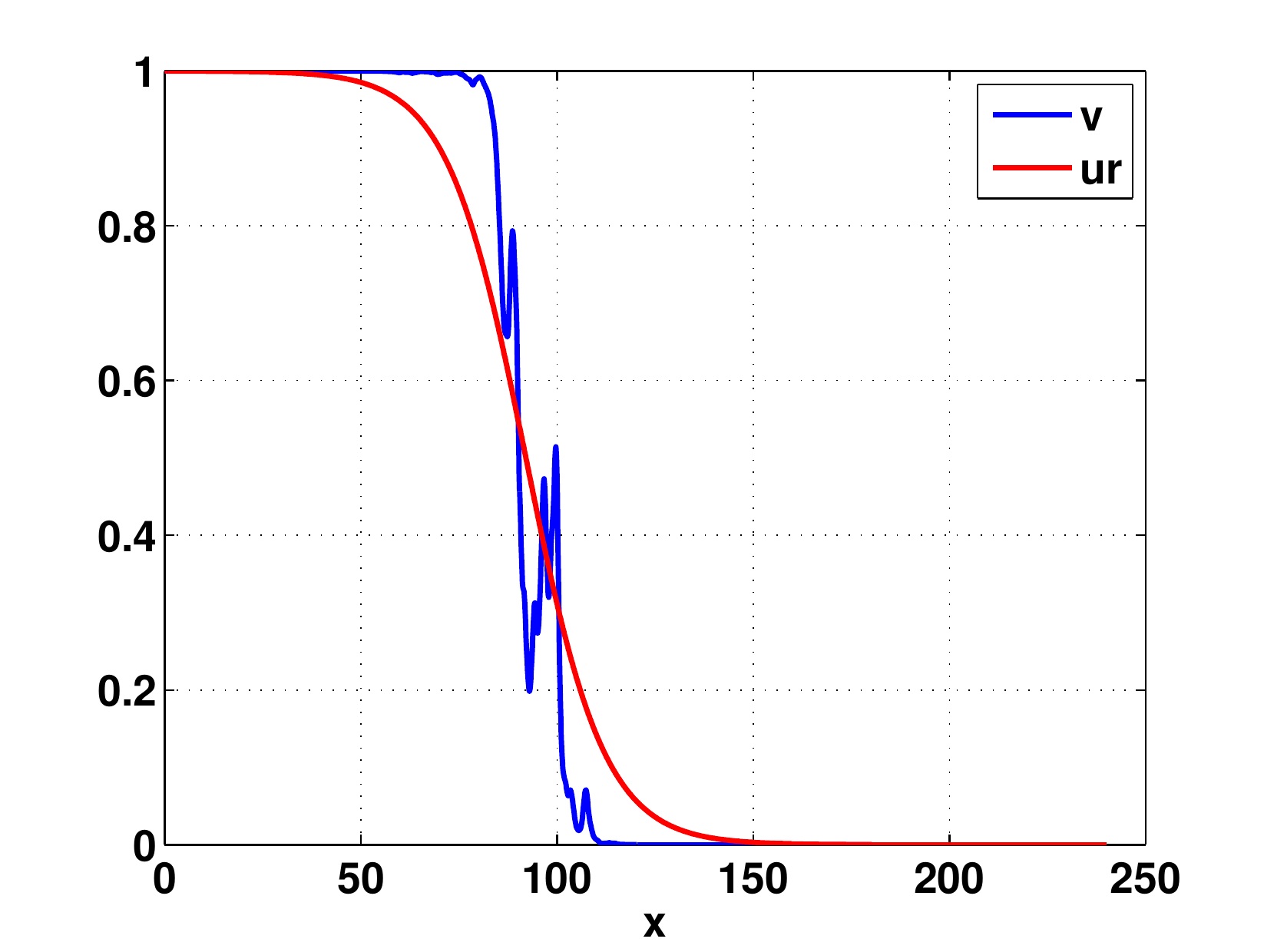}
    \includegraphics*[width=0.4\textwidth]{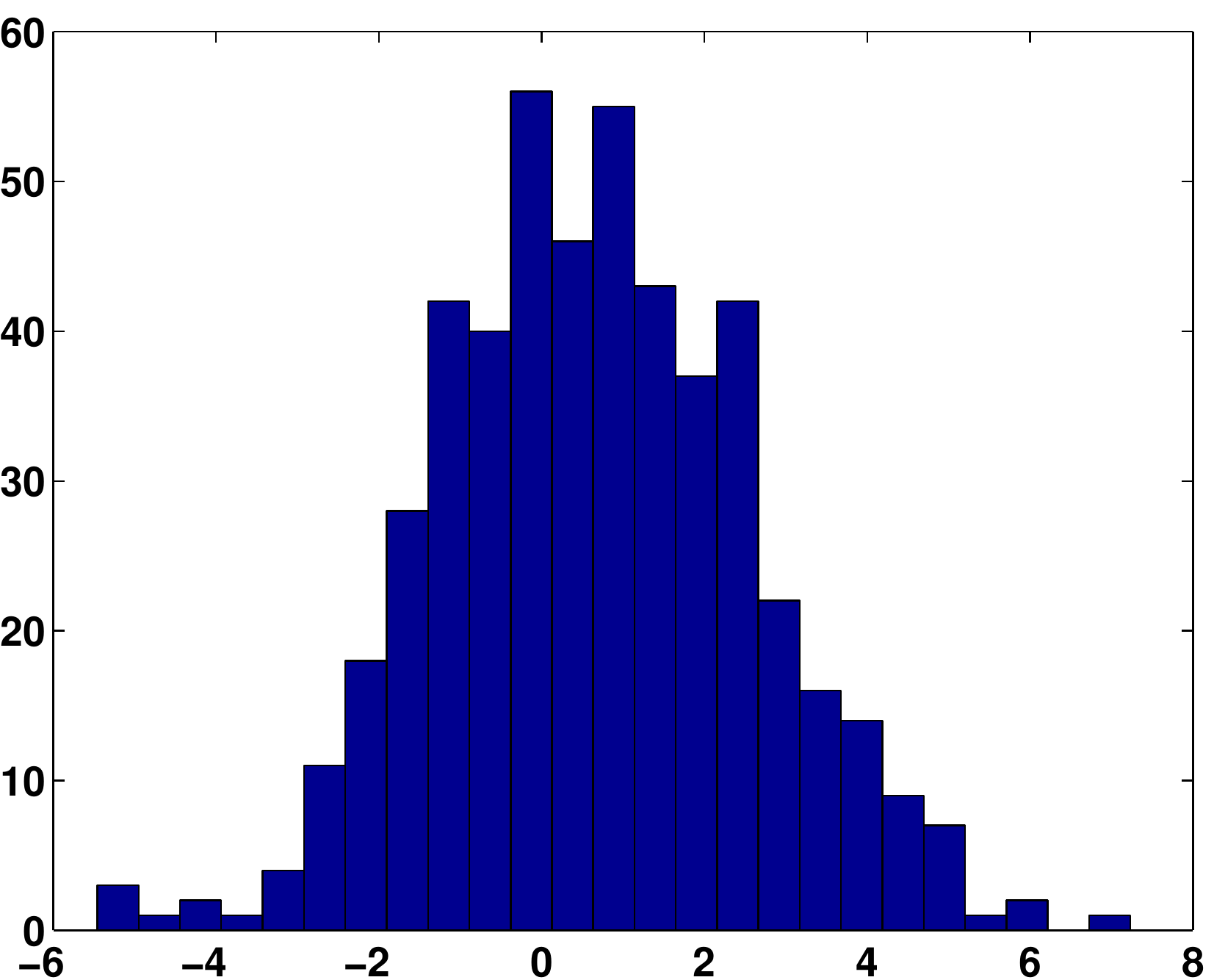}
  \end{center}
  \caption{One realization of the solution with two different reference
    functions $\uh=u_k$ at time $t=50$ (a) and (c) and the
    corresponding different distributions of the instantaneous wave
    speeds $\lambda(t)$ (b) and (d). In
    (a) and (b)  $\kh=1/\sqrt{2}$ and (c) and (d)   $\kh=0.1$. Note
    the smaller variance in (d) with $\kh=1/\sqrt{2}$.}  
  \label{fig:uhfail} 
\end{figure}

To illustrate this we solve the SPDE with  $\alpha=0.25$ and examine
large noise intensity $\mu=1$ combined with a small correlation length
of $\xi=0.5$. 
In \figref{fig:uhfail} (a) with $\uh=u_k$, $k=1/\sqrt{2}$
we see that the width of the computed front is larger than the width
of the reference function $\uh$ and in (b) is plotted the corresponding 
instantaneous wave speed $\lambda(t)$ - with time average $0.6483$
and variance of $6.4568$. 
In (b) the reference function $\uh=u_\kh$ has $\kh=0.1$ and the width
of the reference is larger than the solution. In (d) is plotted the
corresponding  
instantaneous wave speed $\lambda(t)$ - with time average $0.6091$
and variance over time of $3.8006$. 
In (a) the computation of the minimization fails at a later time
($T\approx 55$) however in (b) computation was continued to $T>100$.
In (a) the minimization of the $L^2$ norm 
is dominated by the random
fluctuations in the front which is avoided with a template function
with larger support.
Provided the width of the reference function is comparable or
larger than the width of the front the computations are robust
although convergence rates of the wave speed can be much slower for poor
choices of the reference function. 
%

\subsection{Effects of Stratonovich and \Ito noise}
\label{sec:ItoStrat}
Accurate numerical calculations are notoriously difficult in the
deterministic case when the wave profile depends on the leading
profile of the wave, see for example \cite{ElmrVnVlck} for the Nagumo
equation or \cite{Qiu+Sln:98} for the Fisher equation. 
We consider from now on initial data
that converges to the minimum speed wave in the deterministic case and
take initial data $u_0=u_k$ close to a step function with $k=50$.
We examine the effects on wave speed and support of the front from
changing the noise intensity and correlation length for both
Stratonovich and Ito noise.

\begin{figure}[hbt]
\begin{center}
     (a) \hspace{0.45\textwidth} (b)  \\
  \includegraphics*[width=0.45\textwidth]{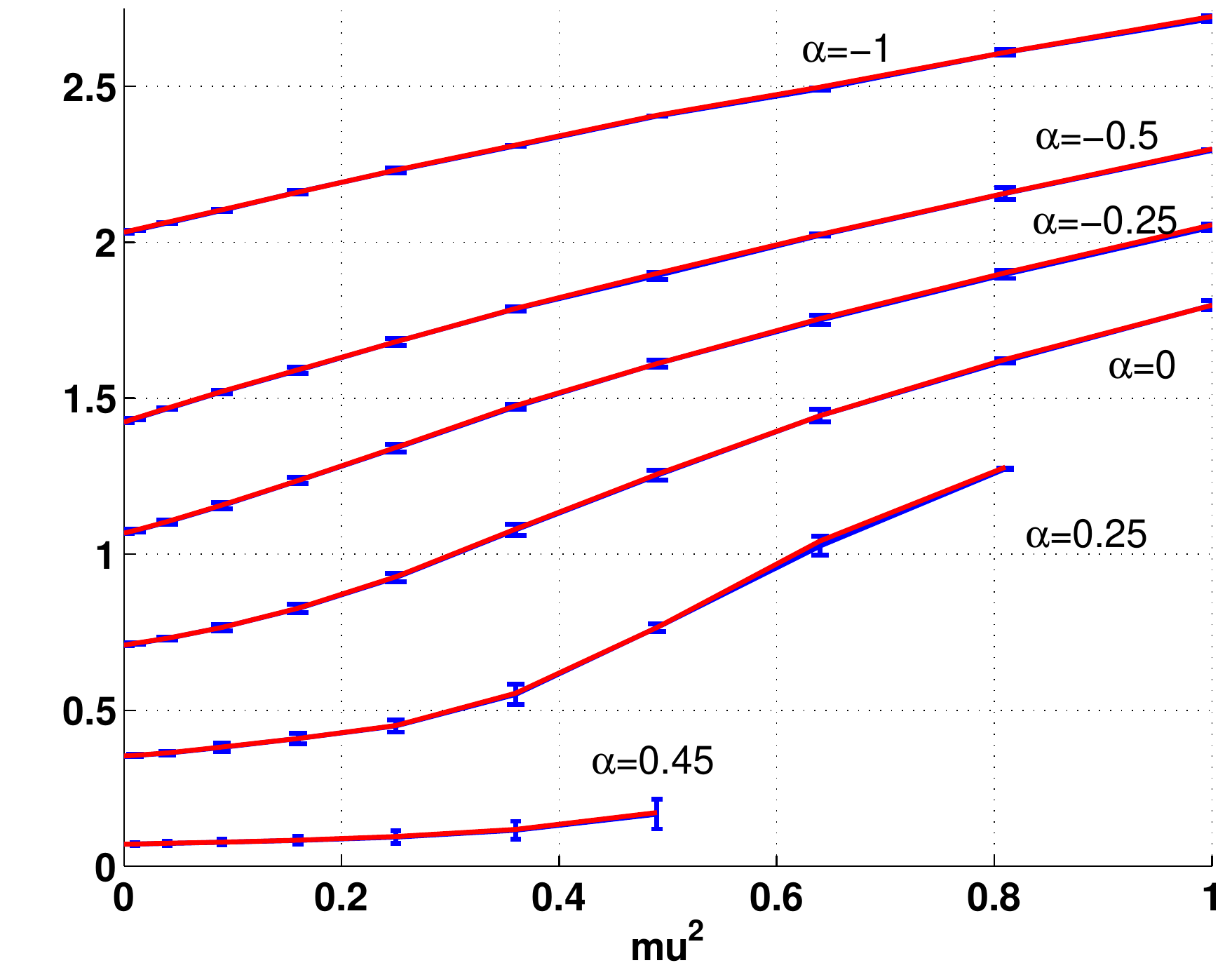}
  \includegraphics*[width=0.45\textwidth]{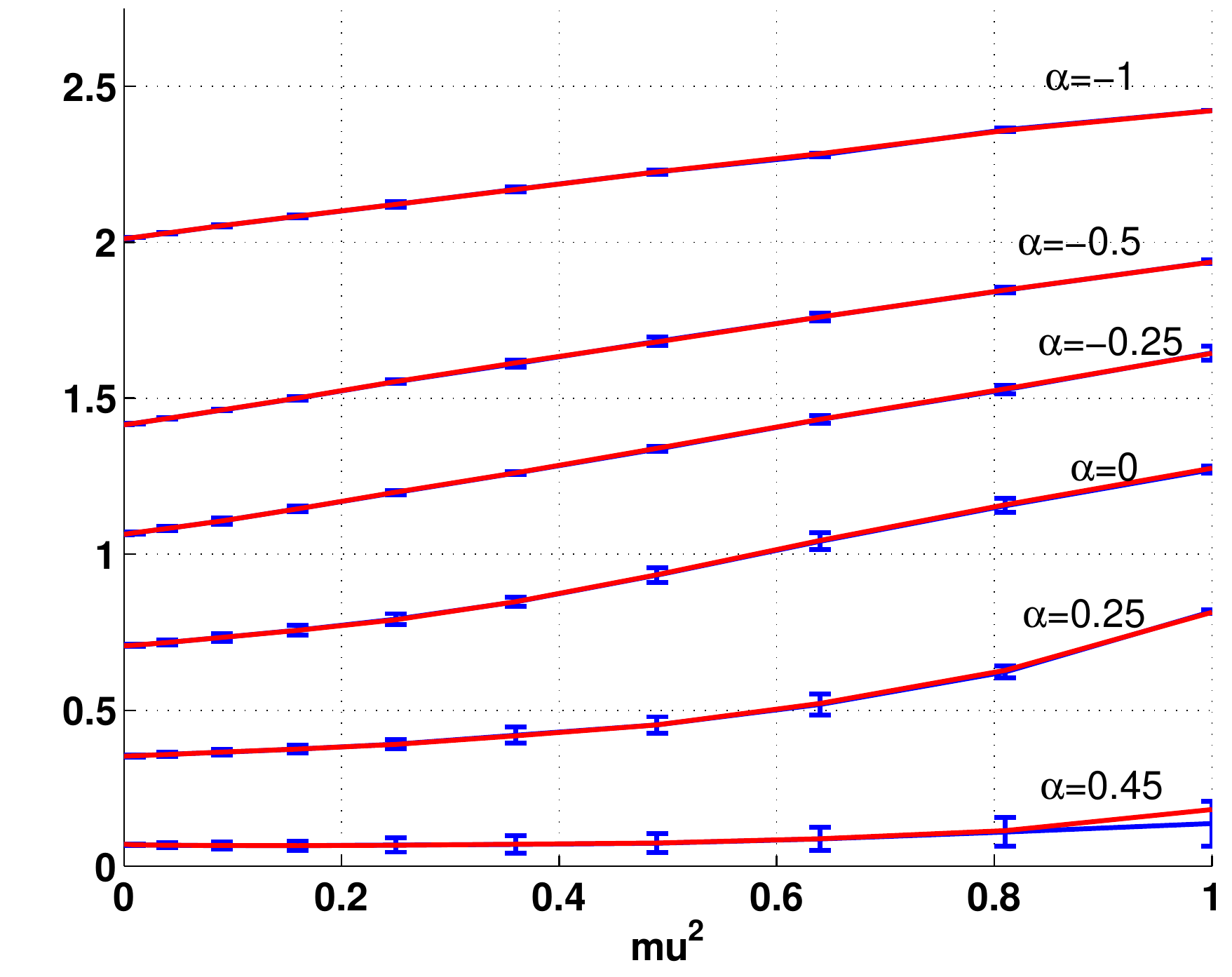}\\
  (c) \hspace{0.45\textwidth} (d)  \\
  \includegraphics*[width=0.45\textwidth]{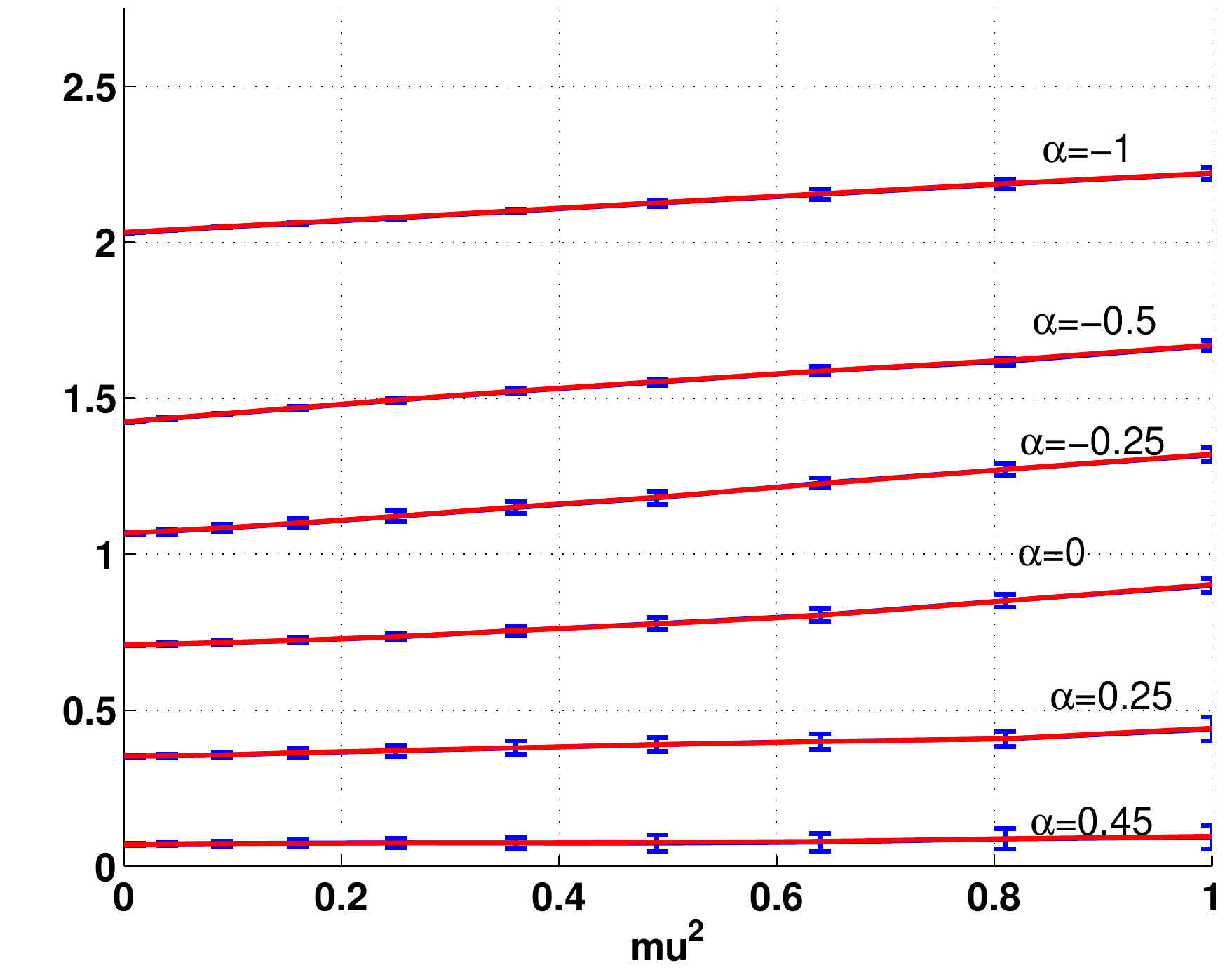}
  \includegraphics*[width=0.45\textwidth]{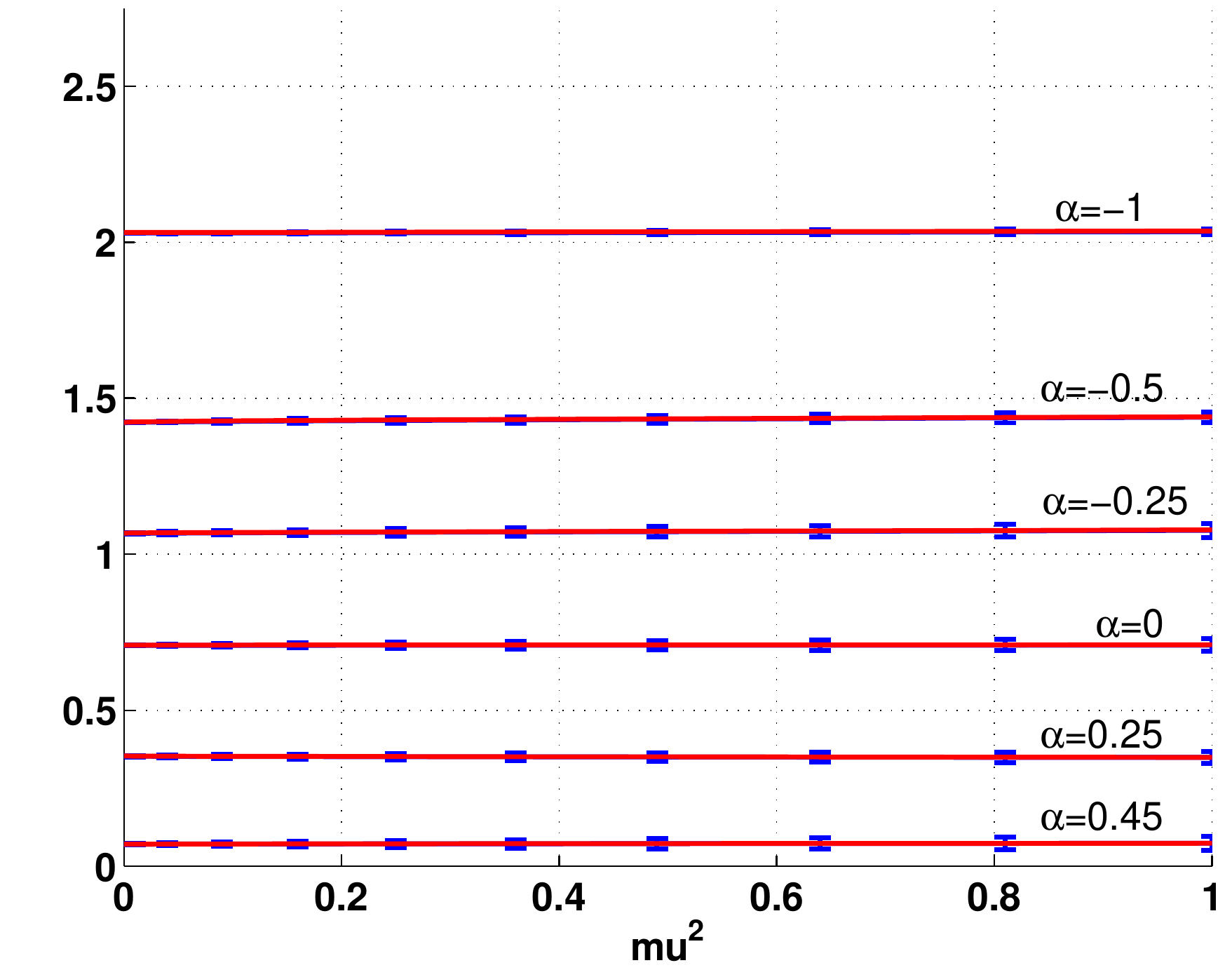}
  \caption{Wave speeds $\LAM$ and $\Lambda_c$ for increasing
    Stratonovich noise intensity and different spatial correlation
    lengths (a) $\xi=0.1$, (b) $\xi=0.5$, (c) $\xi=1$ and (d)
    $\xi=10$. Increasing the noise intensity increases the expected
    wave speed where as increasing the 
    correlation length decreases the expected wave speed.} 
\label{fig:stratlam}
\end{center}
\end{figure}

First we examine the effects of Stratonovich noise on the travelling
wave in the Nagumo equation.
\figref{fig:stratlam} shows wave
speed as noise intensity $\mu$ increases for four different
correlation lengths $\xi=0.1,0.5,1$ and $10$. On each plot are plotted
different nonlinearities $\alpha=0.3,0.25,0,-0.25,-0.3,-0.5,-1$. Each
point on the plot is an average over $100$ realizations and wave
speeds measured both from minimization $\LAM$ and from the level set
$\Lambda_c$.
In \figref{fig:stratlam} we have plotted the corresponding average widths.
We see that increasing the noise intensity increases the wave speed,
where as increasing the correlation length of the noise decreases the
wave speed. The two effects essentially cancel each other in (d) and
we see no overall effect on the noise intensity on the wave speed.
We can also examine the form of the wave profile.
In \figref{fig:stratlam} we have plotted the corresponding average widths of
the wave as noise intensity $\mu$ and correlation length
$\xi$ are changed. For large noise the width of the waves increase and
this effect is again reduced as the correlation length is increased.
For a spatial correlation length $\xi=\Dx=0.1$
we have an approximation of white noise in space, for this case we see
that for $\alpha=0.45$ and $\alpha=0.25$ the width of the wave 
increases and a larger computational domain is required.
%
%
\begin{figure}[hbt]
\begin{center}
     (a) \hspace{0.48  \textwidth} (b)  \\
  \includegraphics*[width=0.45\textwidth]{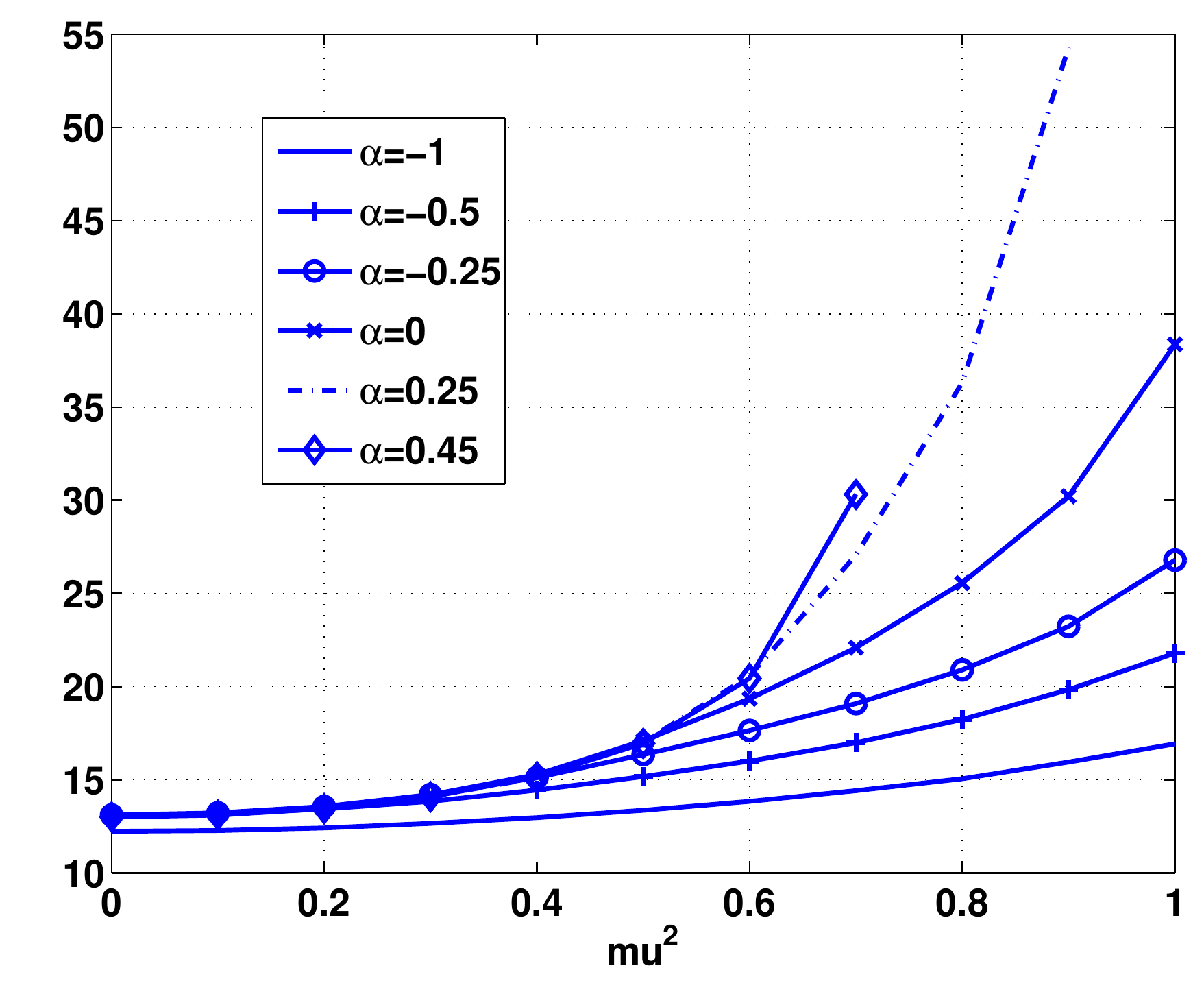}
  \includegraphics*[width=0.45\textwidth]{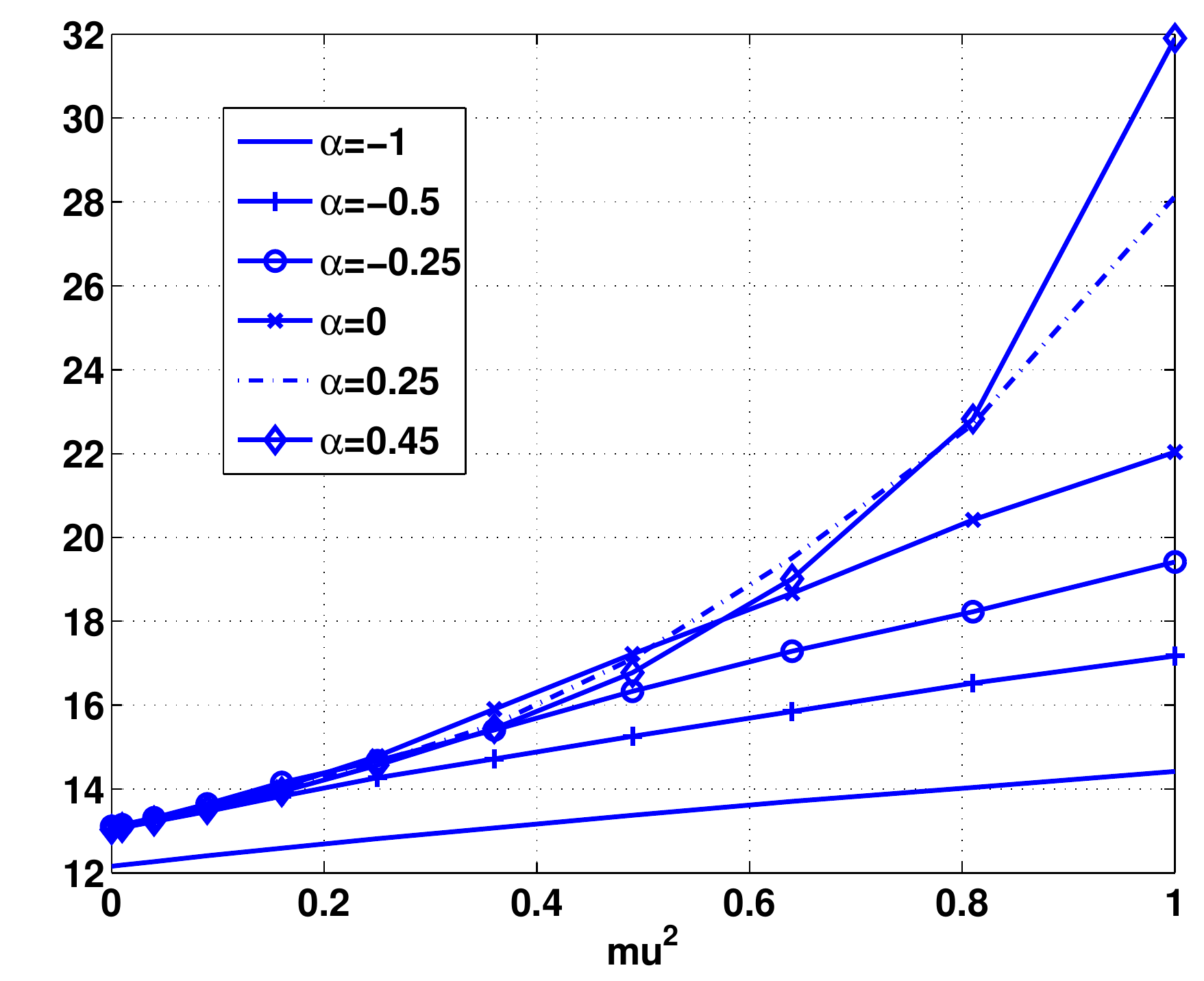}\\
  (c) \hspace{0.48  \textwidth} (d)  \\
  \includegraphics*[width=0.45\textwidth]{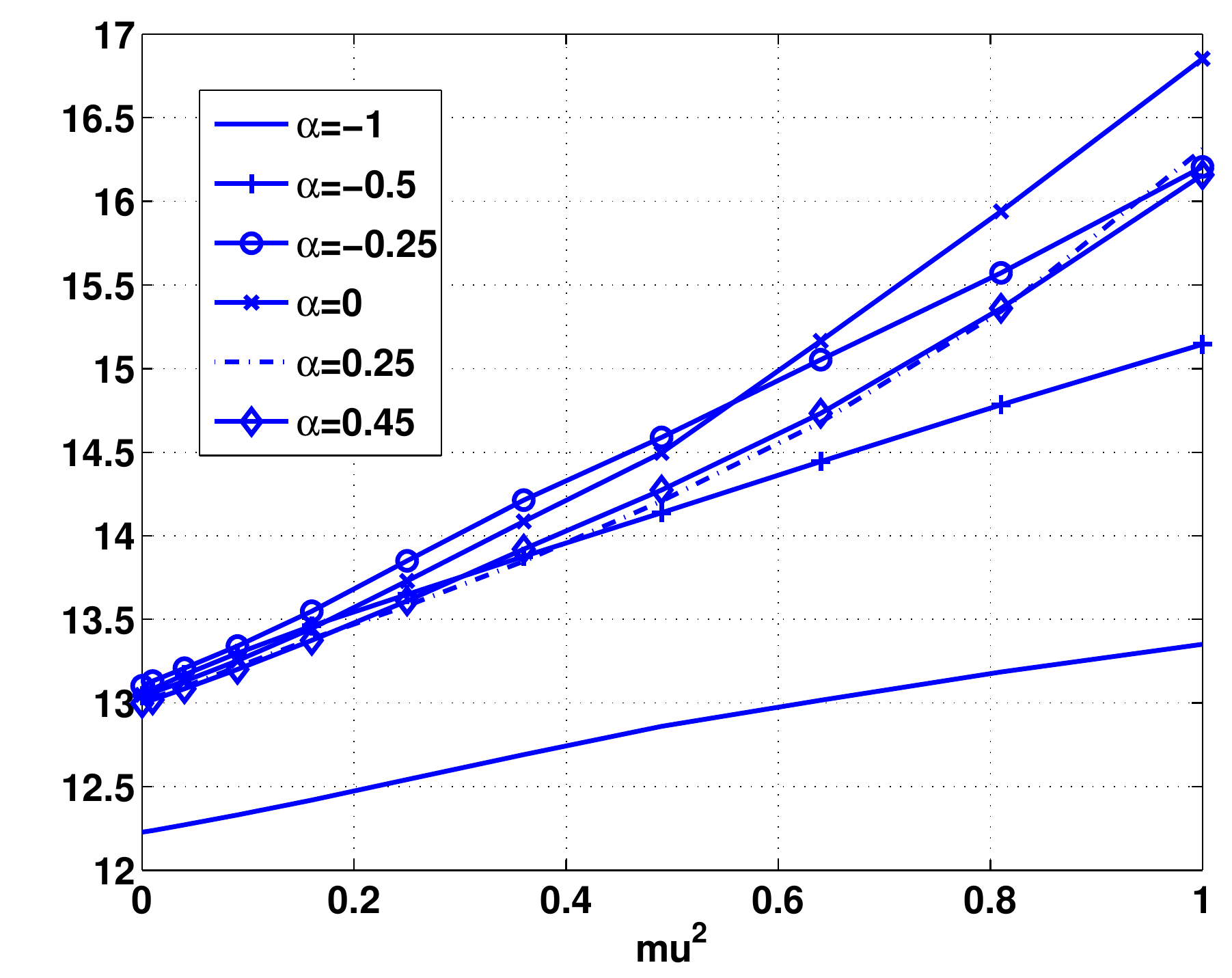}
  \includegraphics*[width=0.45\textwidth]{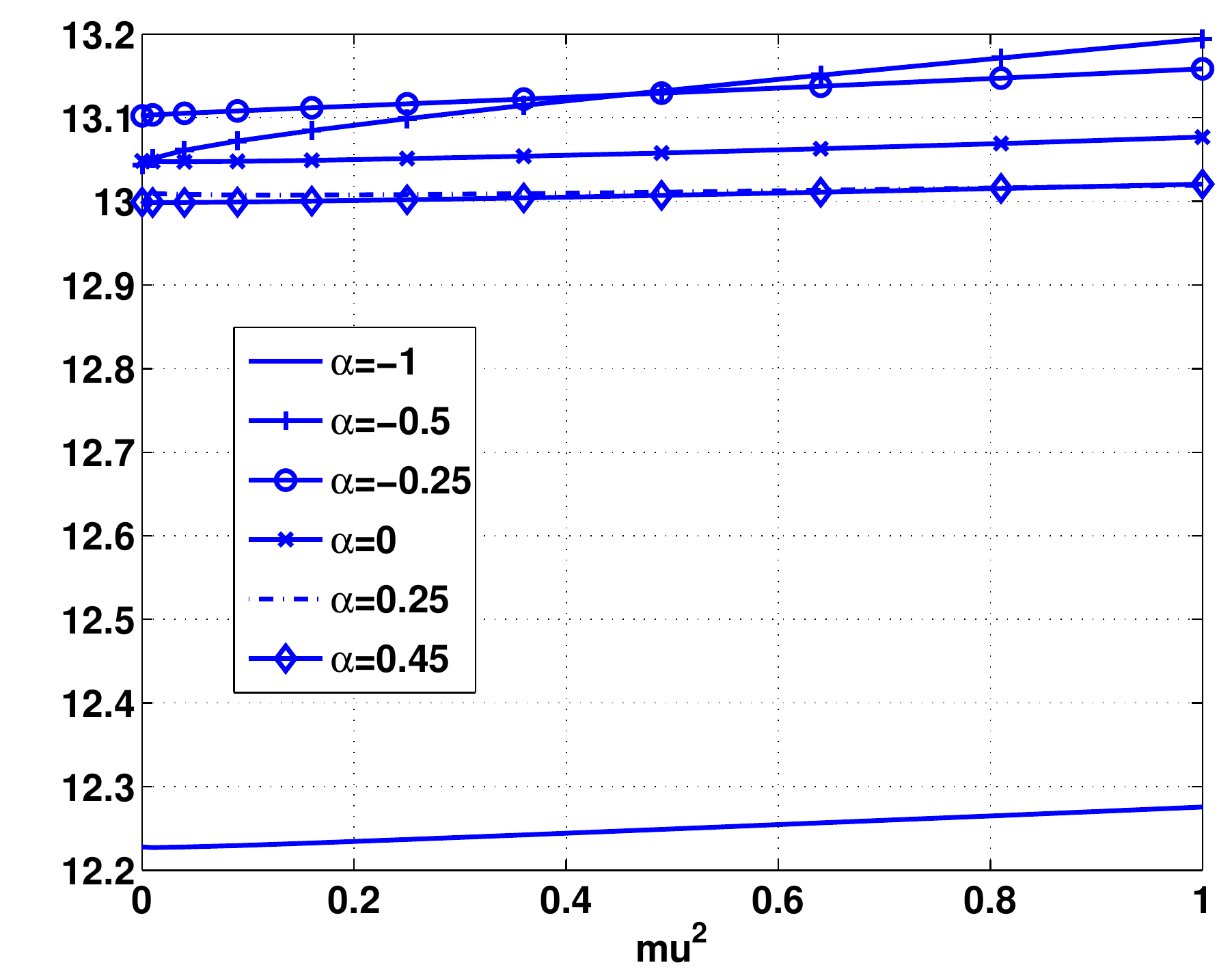}
  \caption{Expected width of the wave for increasing Stratonovich
    noise intensity and  different spatial correlation lengths (a) $\xi=0.1$,
    (b) $\xi=0.5$, (c) $\xi=1$, and (d) $\xi=10$. As the noise
    intensity is increased the expected width of the wave front
    increases where as for fixed intensity increasing the correlation
    length reduces the expected width.}  
\label{fig:stratw}
\end{center}
\end{figure}

\begin{figure}[hbt]
\begin{center}
     (a) \hspace{0.45\textwidth} (b)  \\
  \includegraphics*[width=0.45\textwidth]{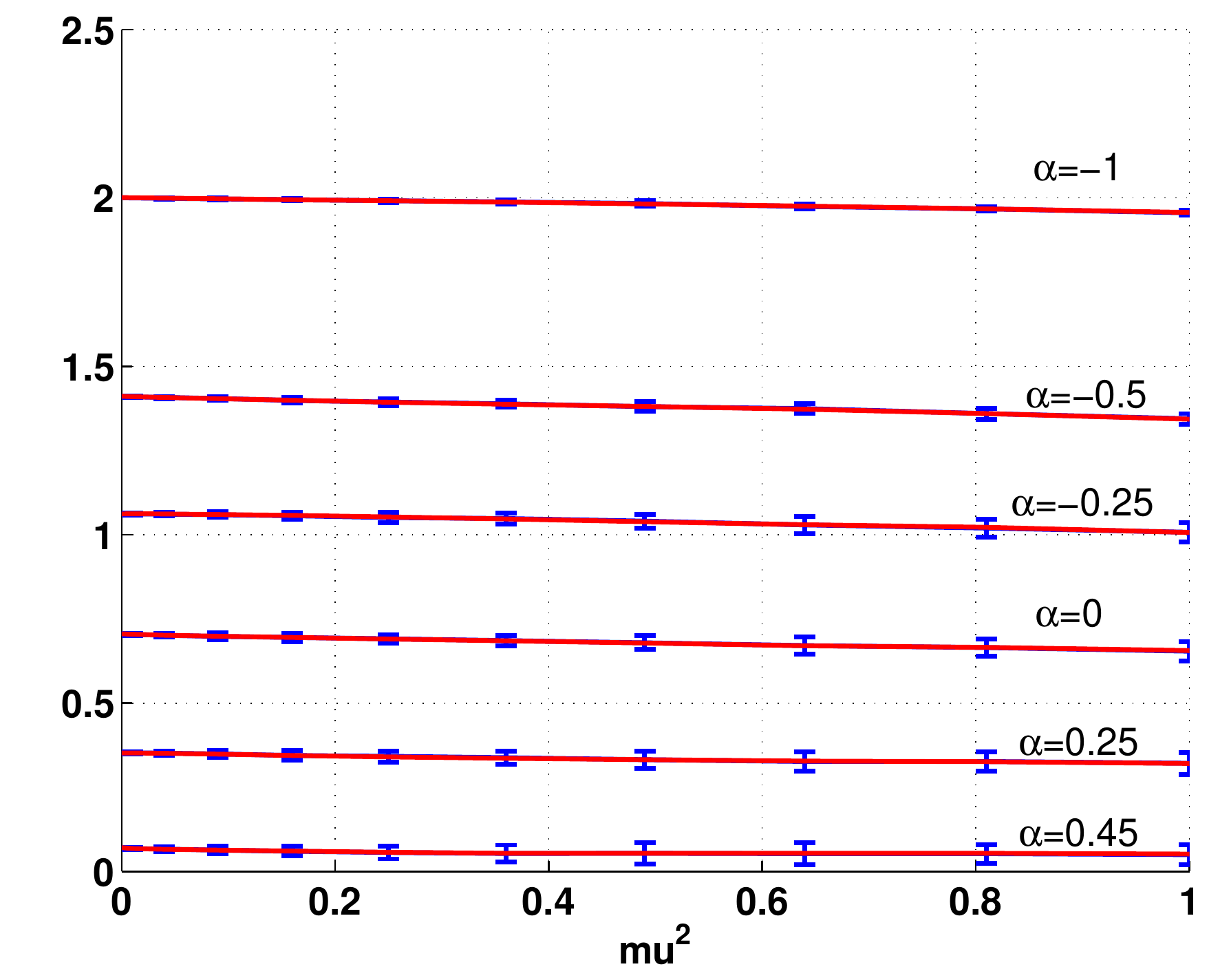}
  \includegraphics*[width=0.45\textwidth]{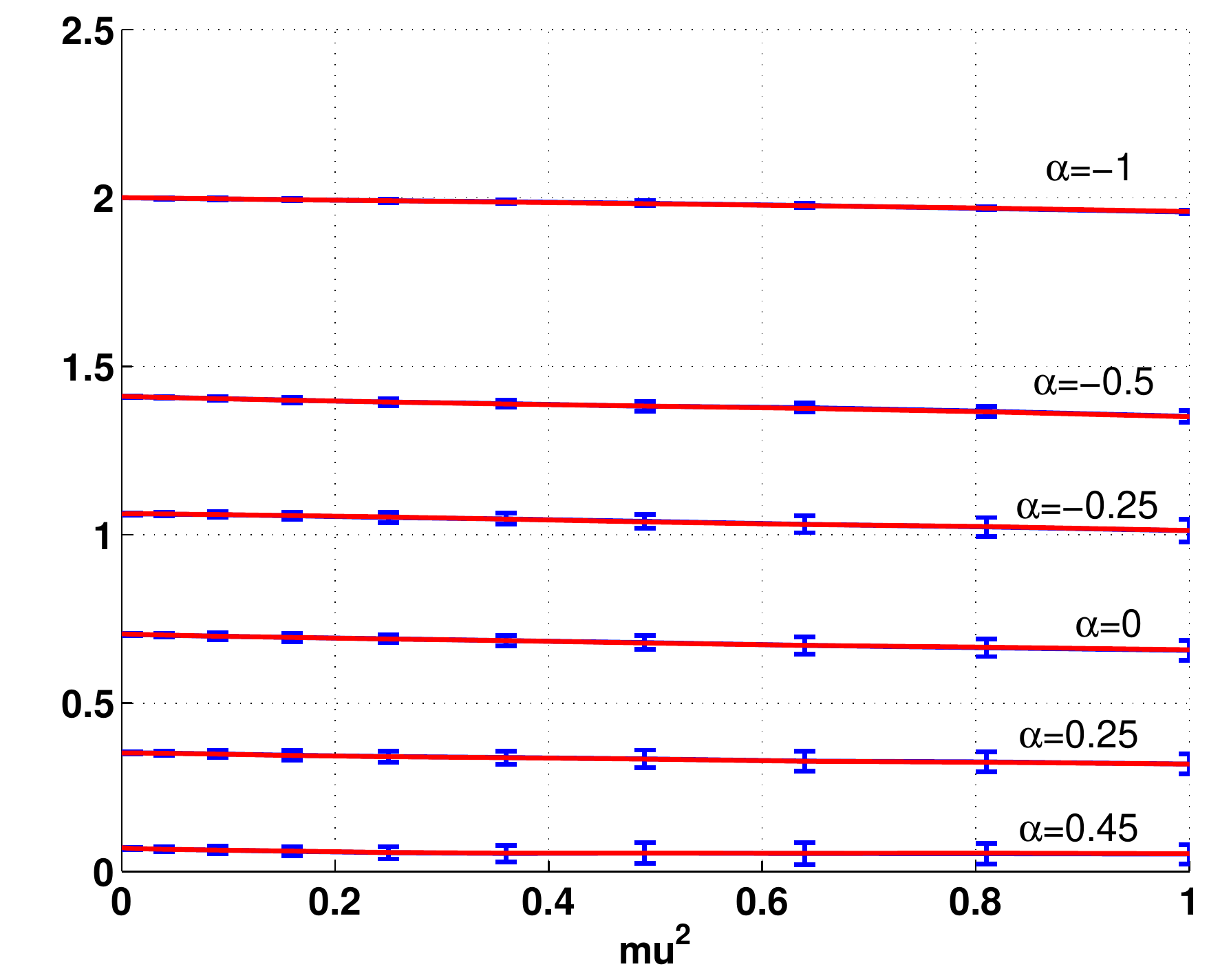}\\
     (c) \hspace{0.45\textwidth} (d)  \\
  \includegraphics*[width=0.45\textwidth]{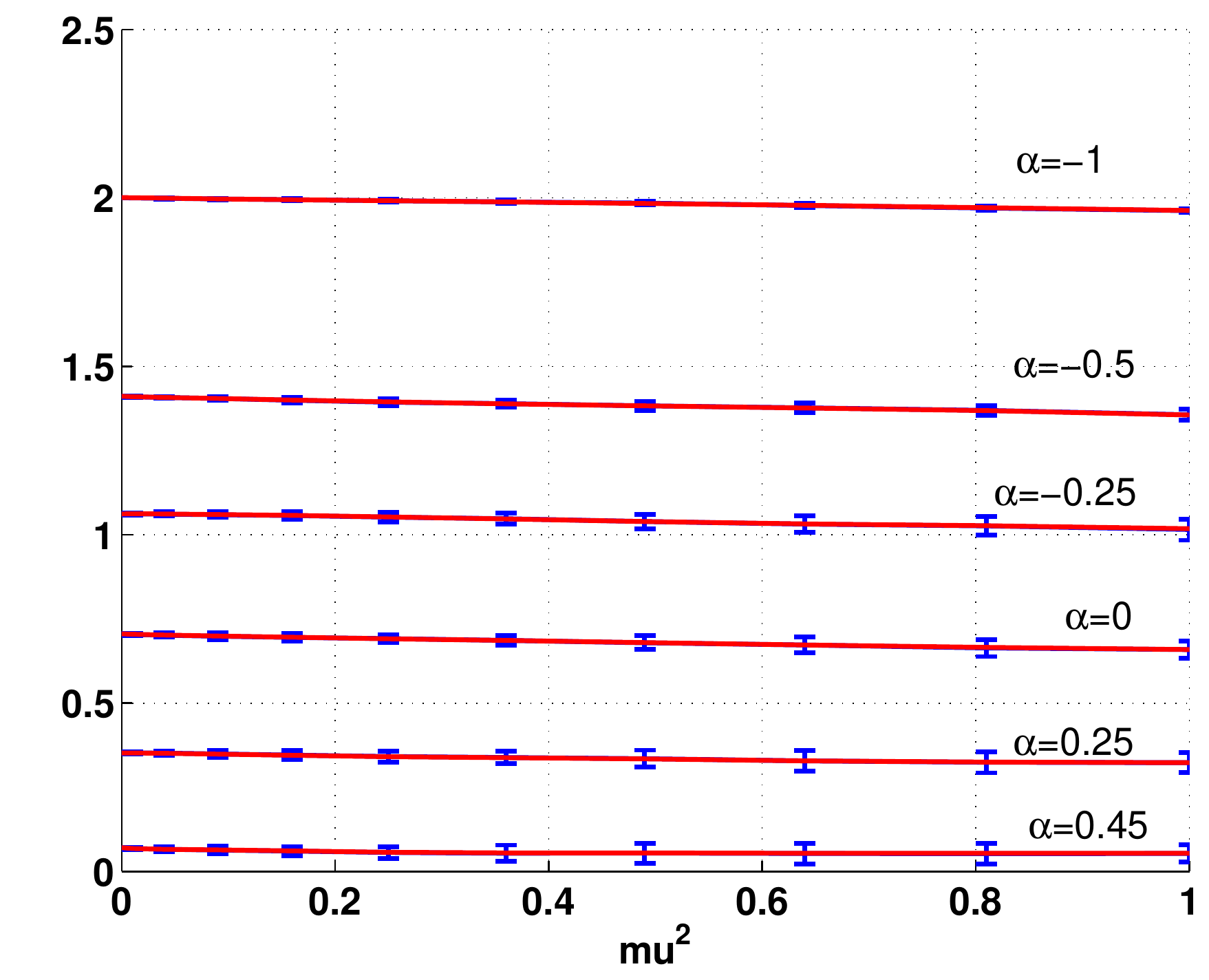}
  \includegraphics*[width=0.45\textwidth]{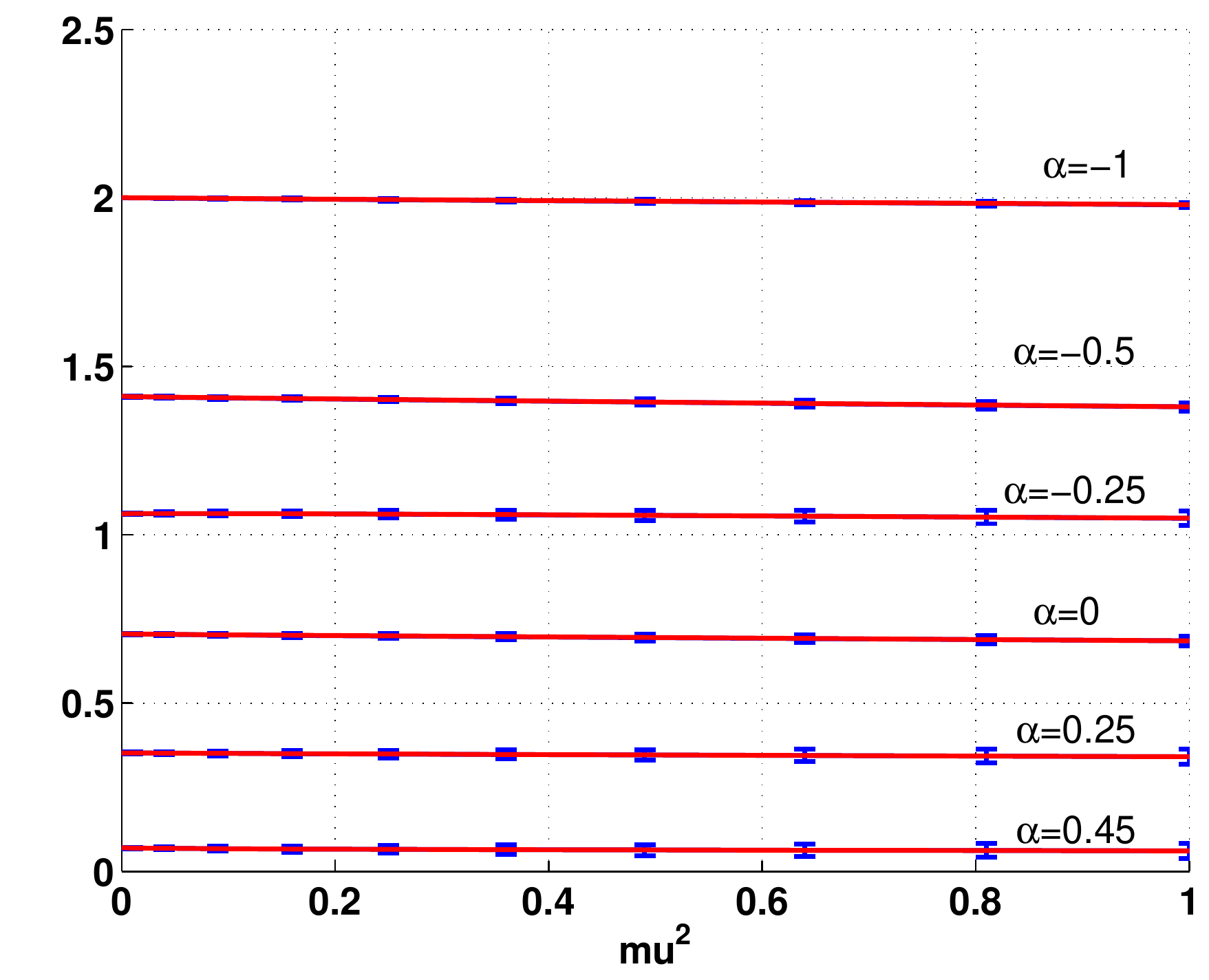}\\
  \caption{Expected wave speeds $\LAM$ and $\Lambda_c$  for increasing
    \Ito noise intensity and 
    different spatial correlation lengths (a) $\xi=0.1$,
    (b) $\xi=0.5$, (c) $\xi=1$ and (d) $\xi=10$. As noise intensity is
  increased we see a slight drop in wave speed and little effect from
  the changing correlation length.} 
\label{fig:itolam}
\end{center}
\end{figure}

\begin{figure}[hbt]
\begin{center}
     (a) \hspace{0.45\textwidth} (b)  \\
  \includegraphics*[width=0.45\textwidth]{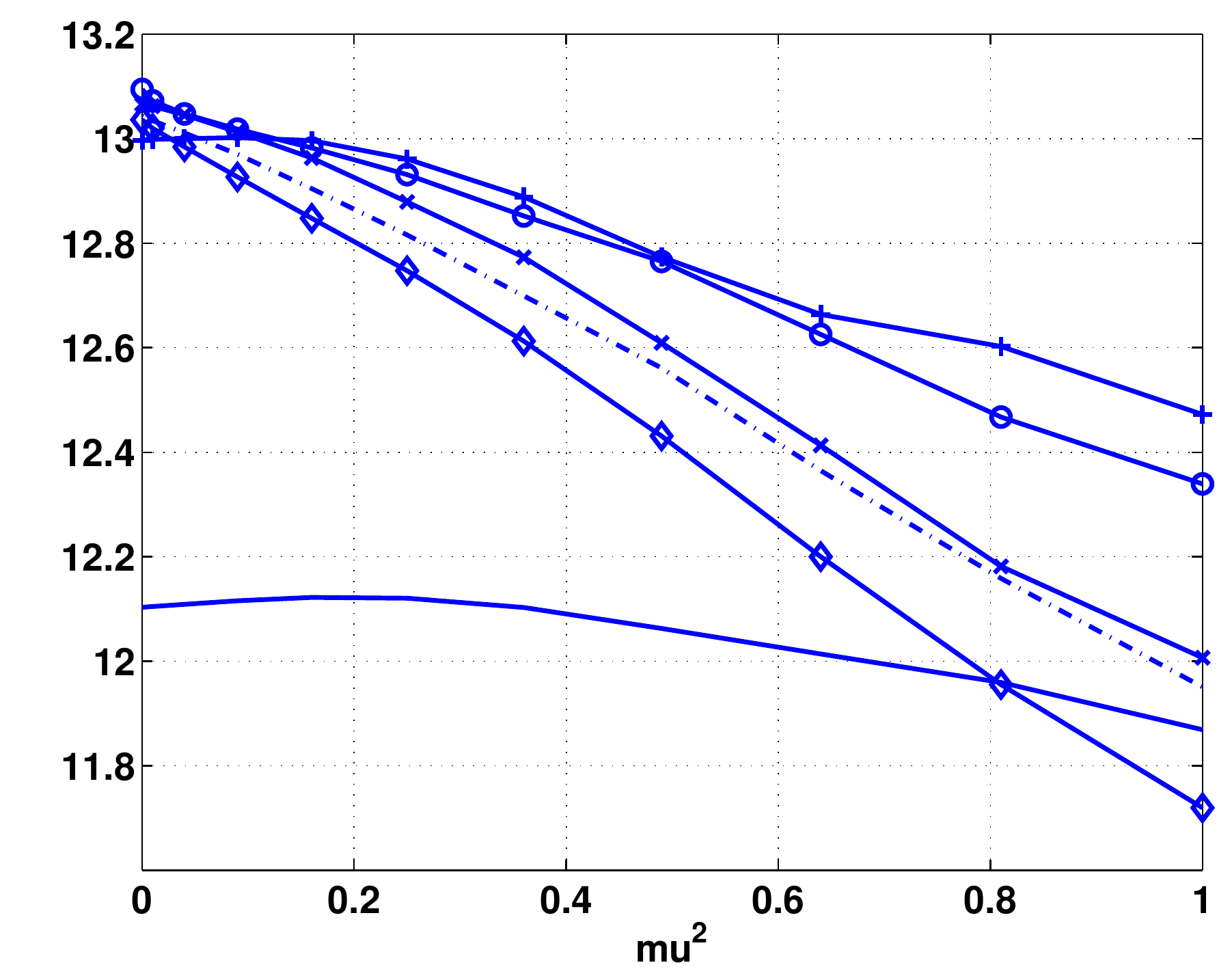}
  \includegraphics*[width=0.45\textwidth]{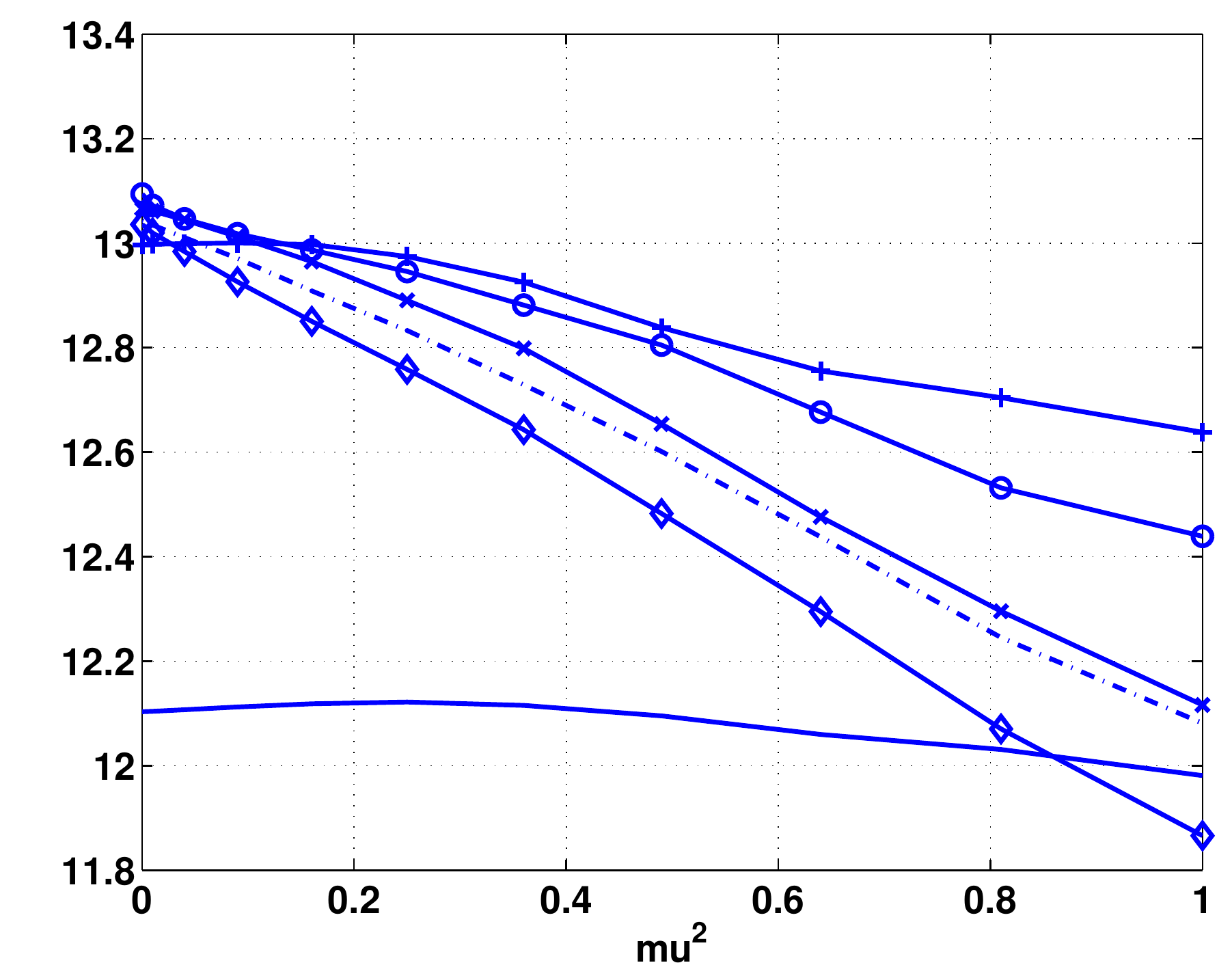}\\
     (c) \hspace{0.45\textwidth} (d)  \\
  \includegraphics*[width=0.45\textwidth]{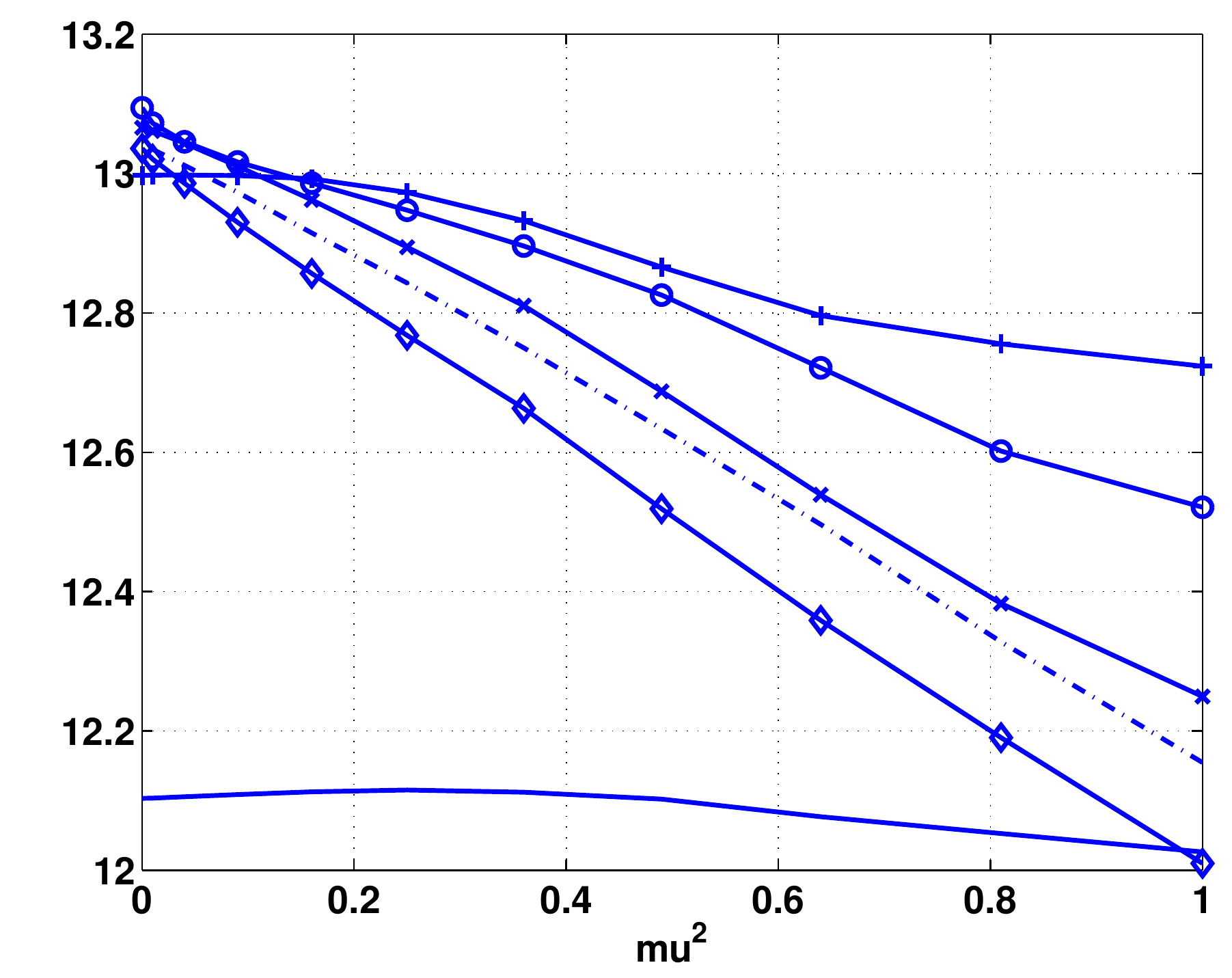}
  \includegraphics*[width=0.45\textwidth]{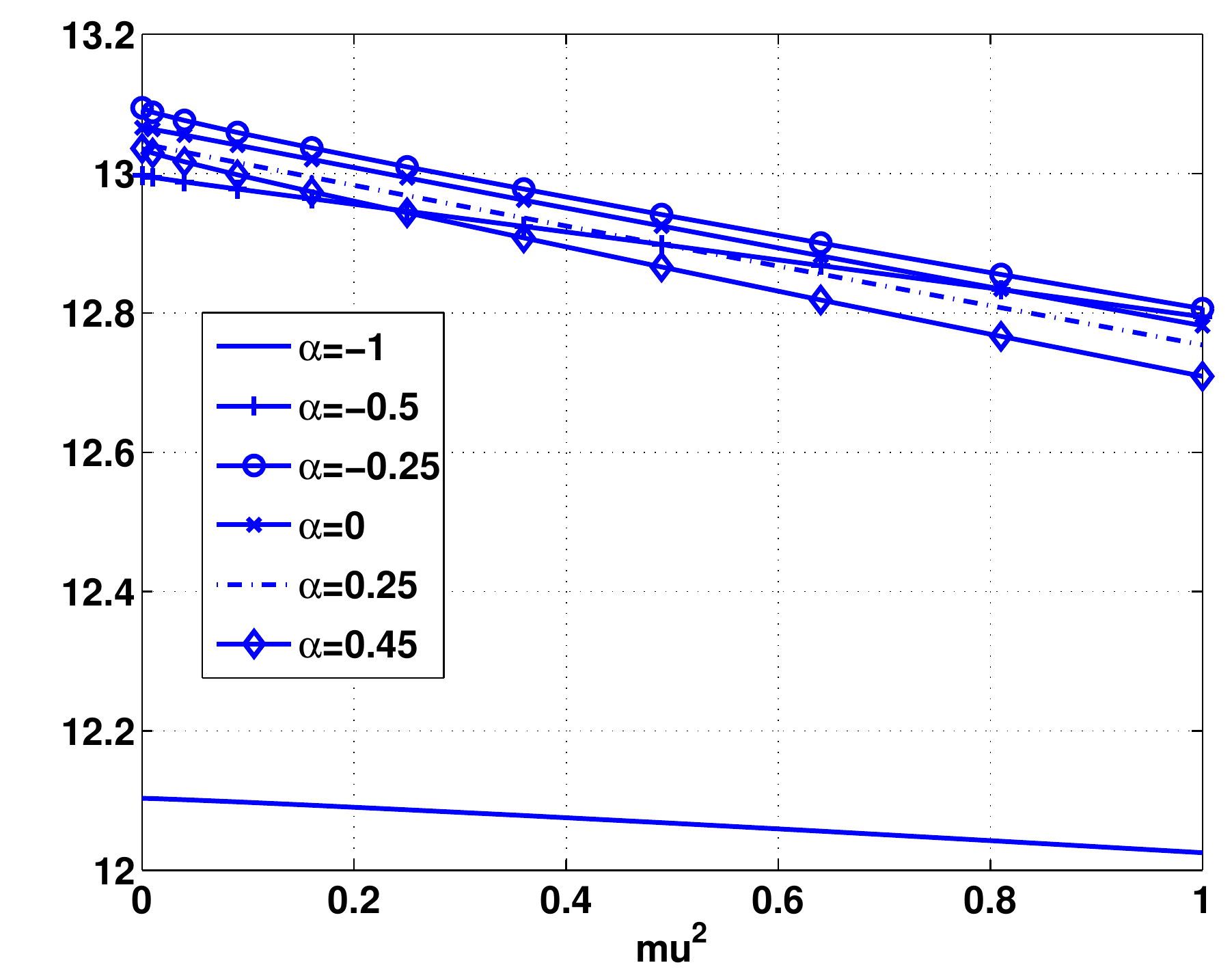}
  \caption{Expected width of waves for increasing \Ito noise intensity and
    different spatial correlation lengths (a) $\xi=0.1$,
    (b) $\xi=0.5$, (c) $\xi=1$ and (d) $\xi=10$. Increasing noise
    intensity narrows the width of the wave and this effect is
    mitigated by increasing the correlation length. Legend for all
    four plots is given in (d).} 
\label{fig:itow}
\end{center}
\end{figure}

For \Ito noise the effect of the noise on wave speed and width of the
waves is less pronounced, see \figref{fig:itolam} for the wave speed
and \figref{fig:itow} for the corresponding width of the waves. 
We see that for large noise, in contrast to the Stratonovich case, a
slight drop in the wave speed for a
correlation length $\xi<10$. As we change the noise intensity we see 
(for most nonlinearities) a drop in the width of the wave -- and so
the front is steeper on average and the effect is more pronounced for
shorter correlation lengths.

\subsection{Computations using averaged quantities}
In general computing wave profiles using averaged quantities leads to the wave
being 'polluted' by the spread of the individual waves (see
\cite{GrciaOjlvo+Sncho}).

\begin{figure}[hbt]
  \begin{center}
    (a) \hspace{0.45\textwidth} (b)  \\
    \includegraphics*[width=0.42\textwidth]{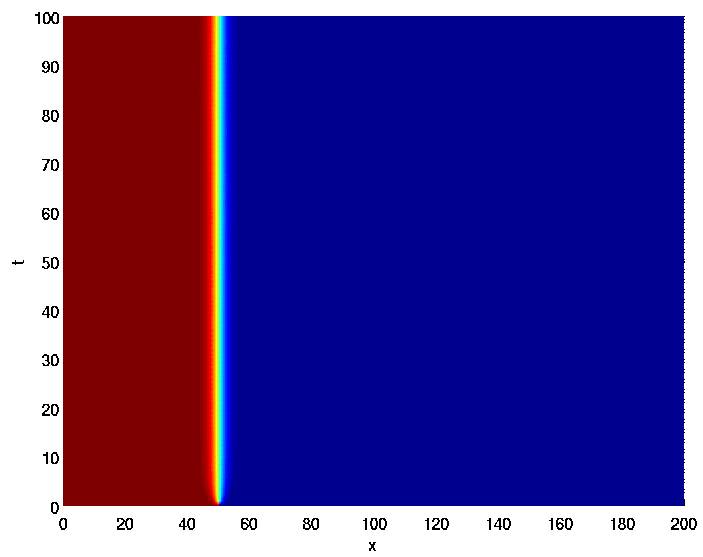}
    \includegraphics*[width=0.42\textwidth]{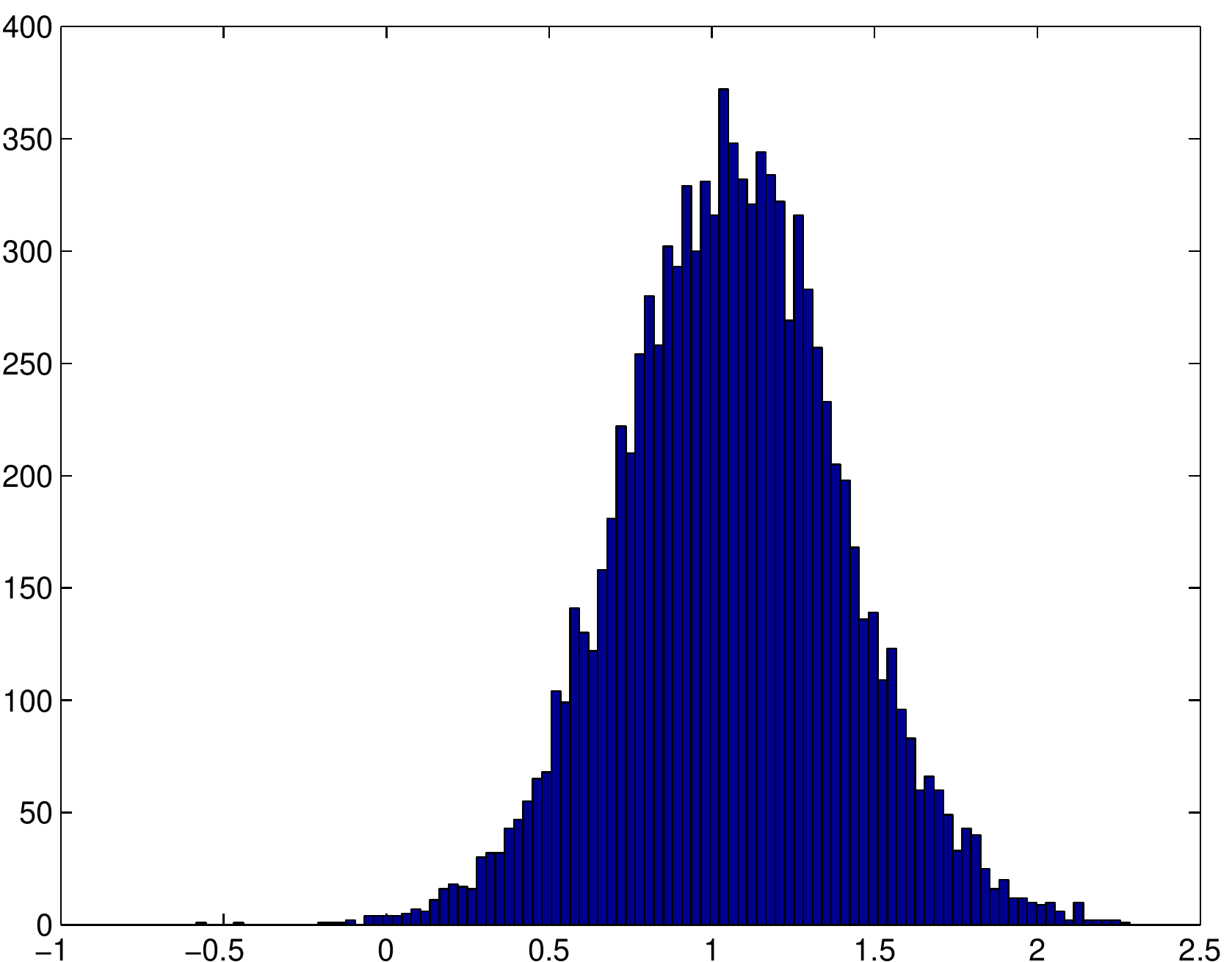} \\
    (c) \hspace{0.45\textwidth} (d)  \\
    \includegraphics*[width=0.42\textwidth]{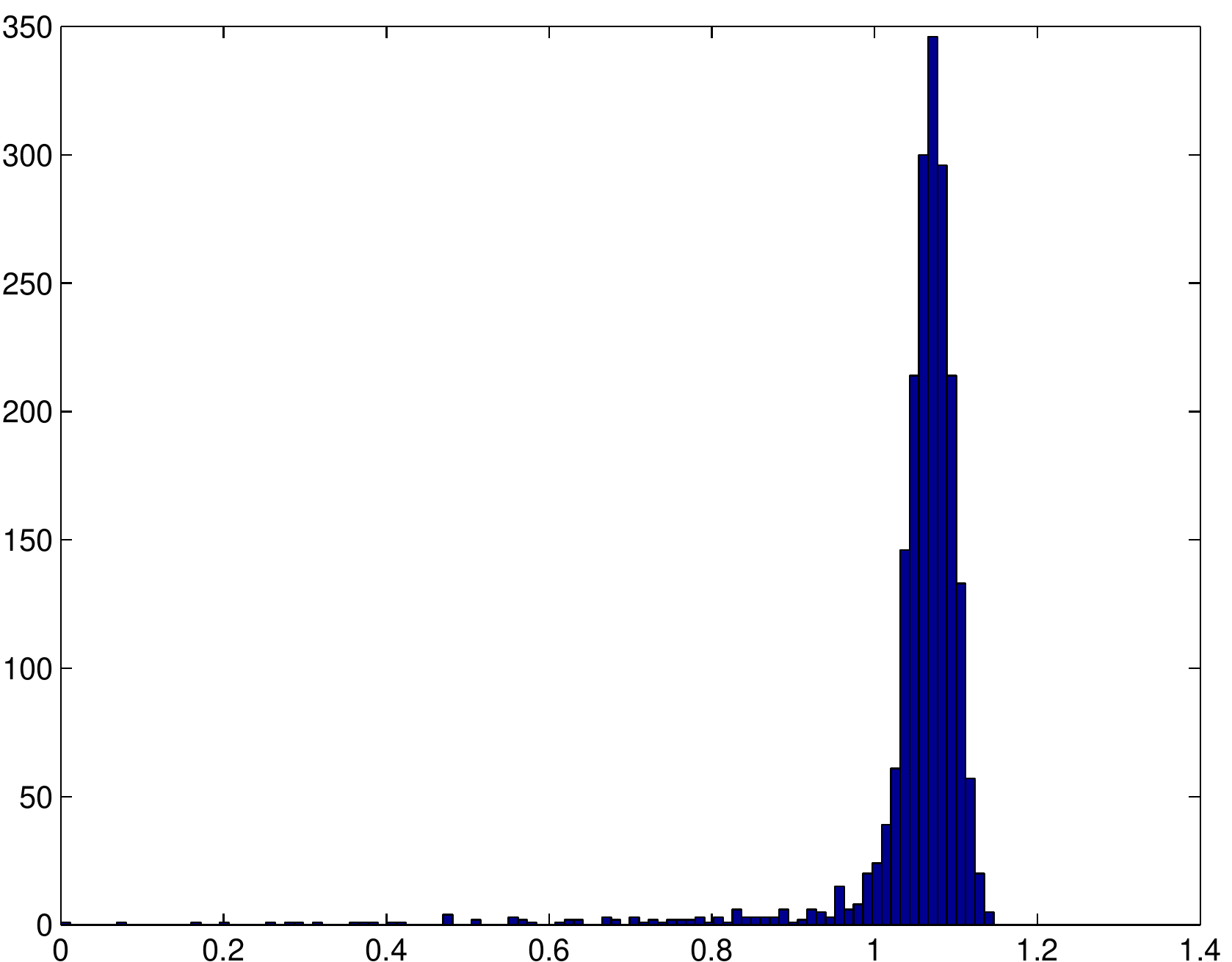}
    \includegraphics*[width=0.42\textwidth]{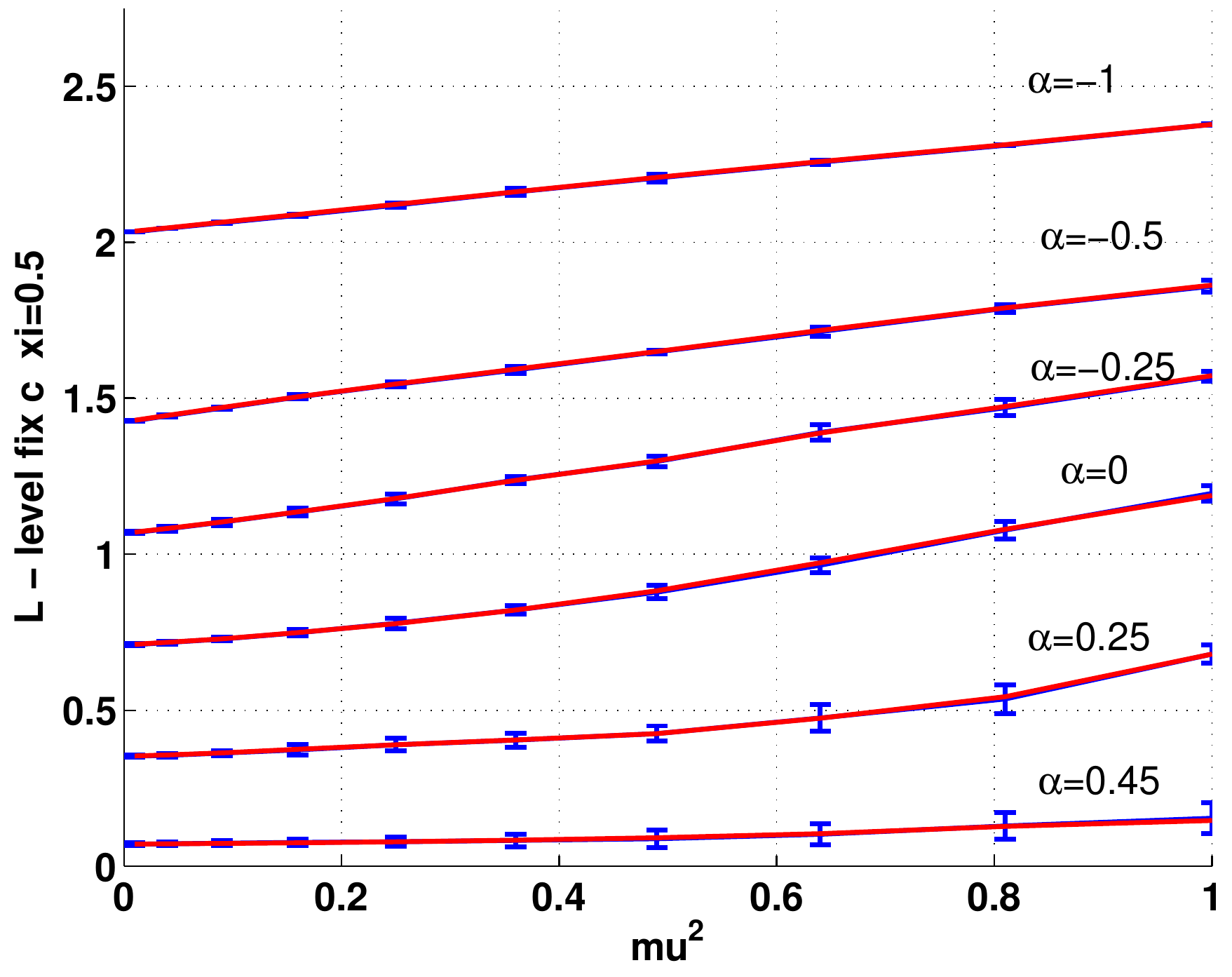} 
  \end{center}
\caption{(a) Mean solution of SPDAE \eqref{eq:spdae.E}. Note that we
  do not observe that the wave front has be spread taking the mean
  instantaneous speed. In (b) is plotted the distribution of $\lambda$
used to freeze the wave in (a). In (c) we plot the distribution from
solving the SPDE using an average wave speed for the minimization and
for the SPDE and in (d) we compare wave speeds over a range on
nonlinearities and noise intensities for $\xi=0.5$.
}
\label{fig:Elm}
\end{figure}

In \eqref{eq:spdae.E} we propose using an expected value of the
instantaneous wave speeds for the SPDAE. We fix a spatial correlation
of $\xi=0.5$. For the SPDAE if we solve with
$\mu=0.1$, $\uh=u_k$ with $k=0.1$ and $u^0=u_{k_0}$, $k_0=1/\sqrt{2}$
and $100$ realizations then we obtain an estimate of a wave speed of
$1.086$. This compares with $\LAM=1.086$ and $\Lambda_c=1.084$ from
solving the SPDE \eqref{eq:spdeStrat}.
If we examine the computed mean solution front we do not observe
spreading of the wave front (see \figref{fig:Elm} (a)). We also note
from (b) that the distribution of $\lambda(t)$ has smaller variance
than that from solving \eqref{eq:stratspdae}.
For the SPDE we can implement a version \eqref{eq:spdae.E} where we
move the reference function using the expected values of the
instantaneous wave speeds. In \figref{fig:Elm} (c) we plot the
distribution of $\lambda$ for same parameters as in (b). The mean
values agree although the distributions are different.
In (d) we see that computed wave speeds using the average
instantaneous speed and wave speeds $\Lambda_c$ computed using the
level set approach are the same over a range of nonlinearities. These
are the same as those computed using  the SPDE, compare to
\figref{fig:stratlam} (b).

\subsection{Additive noise}
\label{sec:add}
We briefly consider the case of additive noise in 
the \SPDE for which, unless
the noise has some special properties, a solution will in general
cease to exist at some finite time. 

\begin{figure}[hbt] 
\begin{center}
  (a) \hspace{0.32\textwidth} (b)  \\
  \includegraphics*[width=0.45\textwidth]{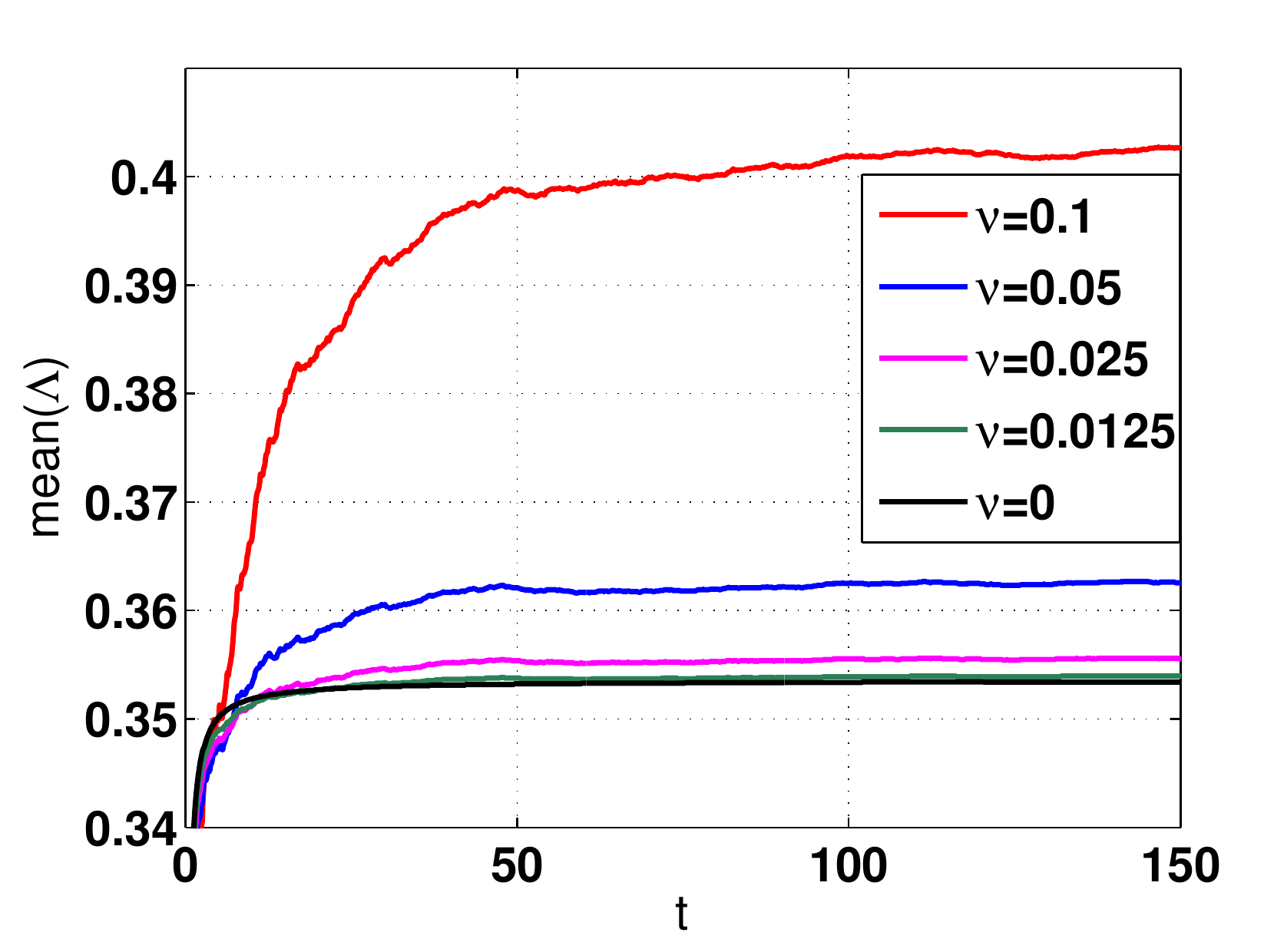}
  \includegraphics*[width=0.45\textwidth]{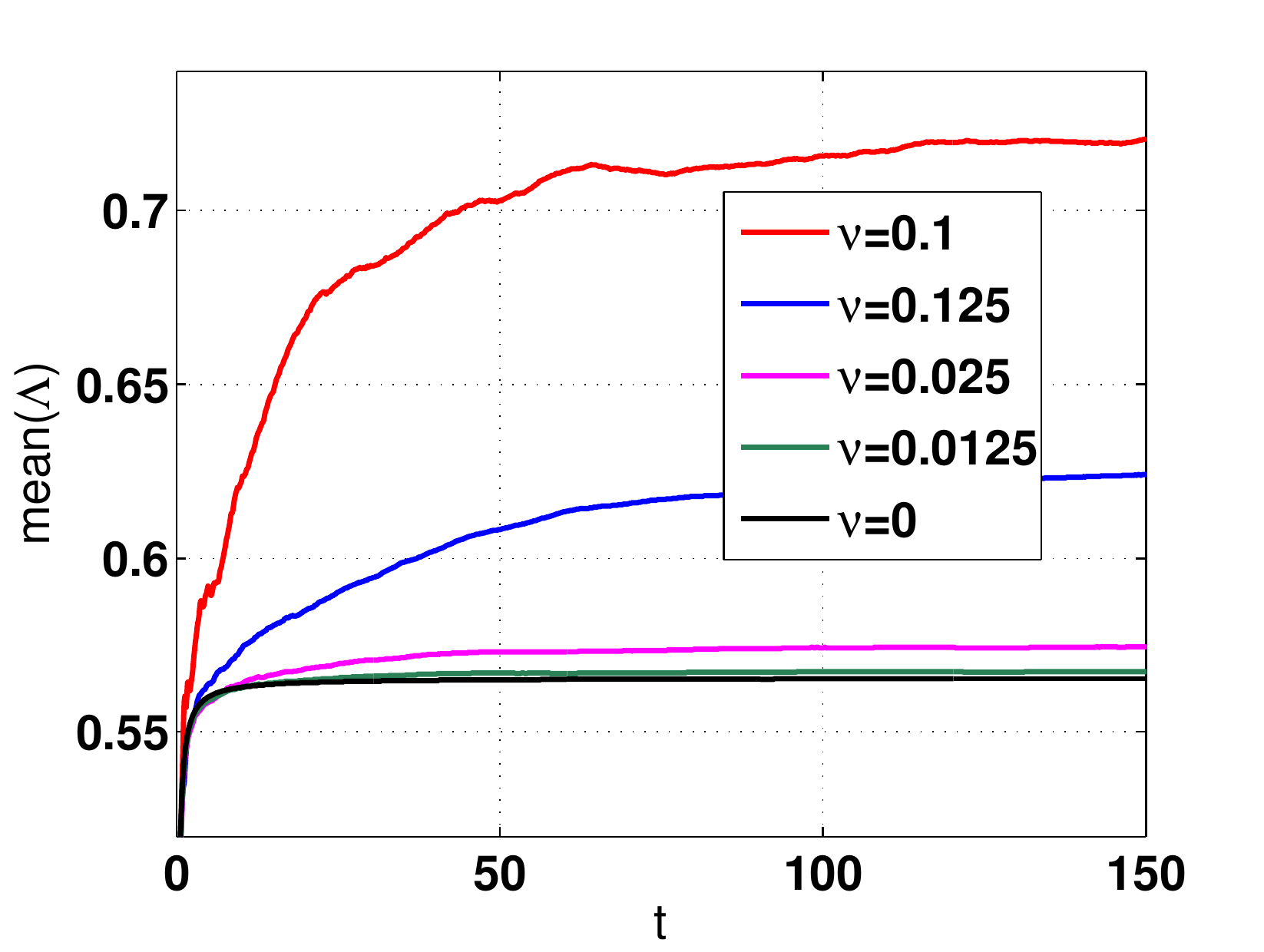}
\caption{Influence of additive white noise on the wave speed $\LAM$ for (a)
  $\alpha=0.25$ and (b) $\alpha=0.1$. In both cases there is a clear
  increase in the wave speed from the deterministic case ($\nu=0$) as
  the noise intensity is increased, with the final wave speeds in
  order of the indicated niose intensities.}
\end{center}
\end{figure}

\begin{figure}[hbt]
\begin{center}
  (a) \hspace{0.45\textwidth} (b)  \\
  \includegraphics*[width=0.45\textwidth]{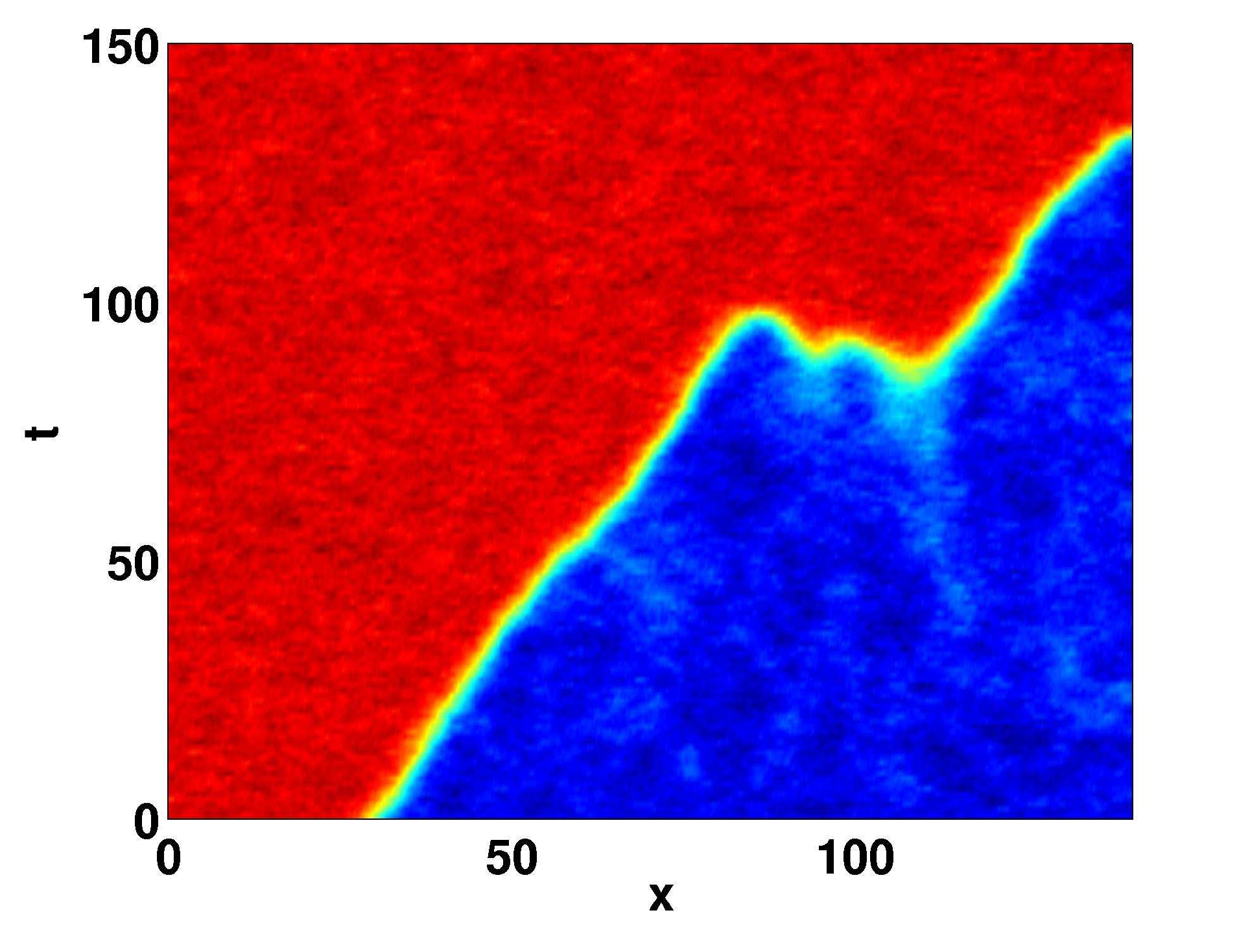}
  \includegraphics*[width=0.45\textwidth]{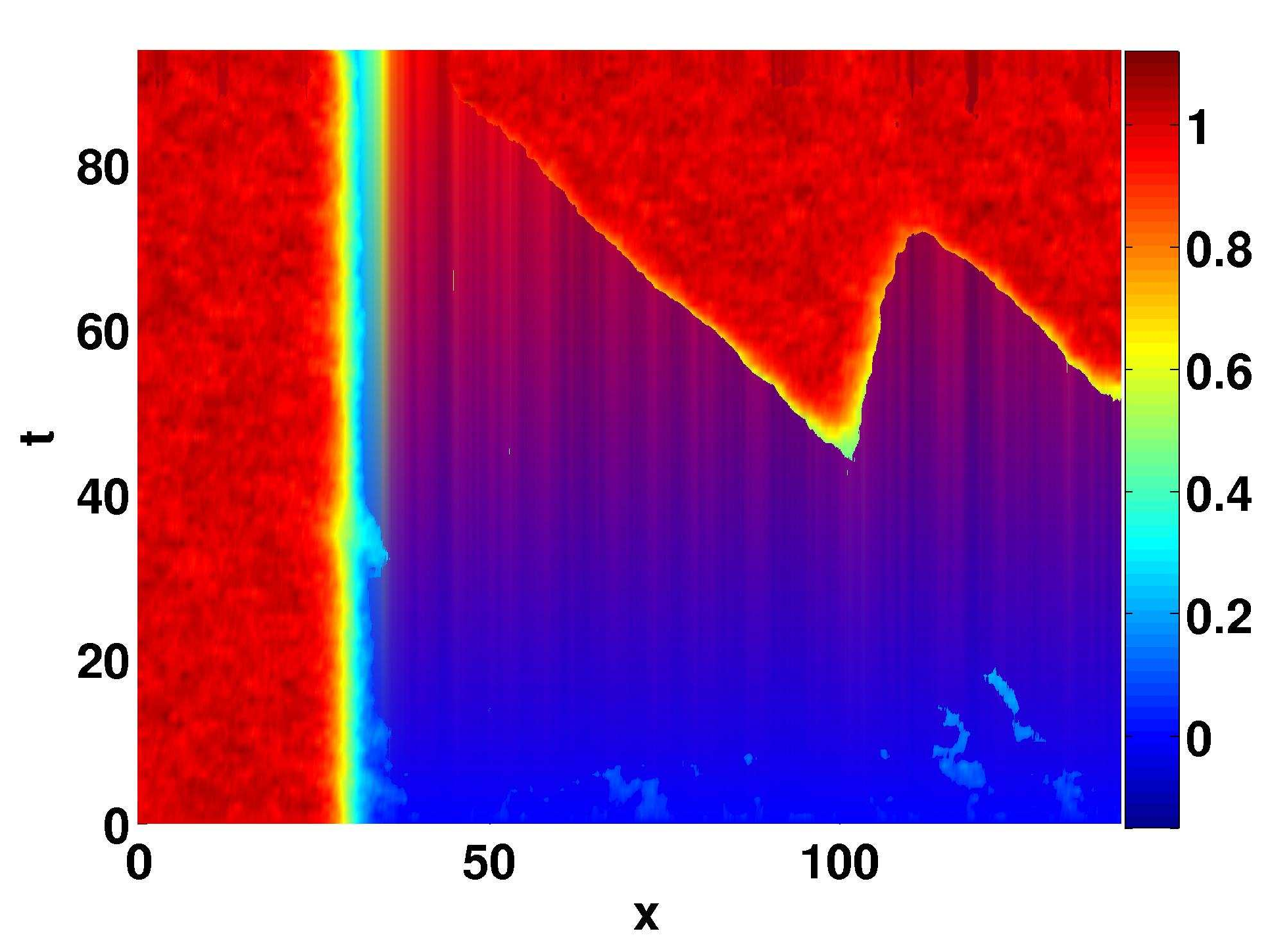}\\
  (c) \hspace{0.45\textwidth} (d)  \\
  \includegraphics*[width=0.45\textwidth]{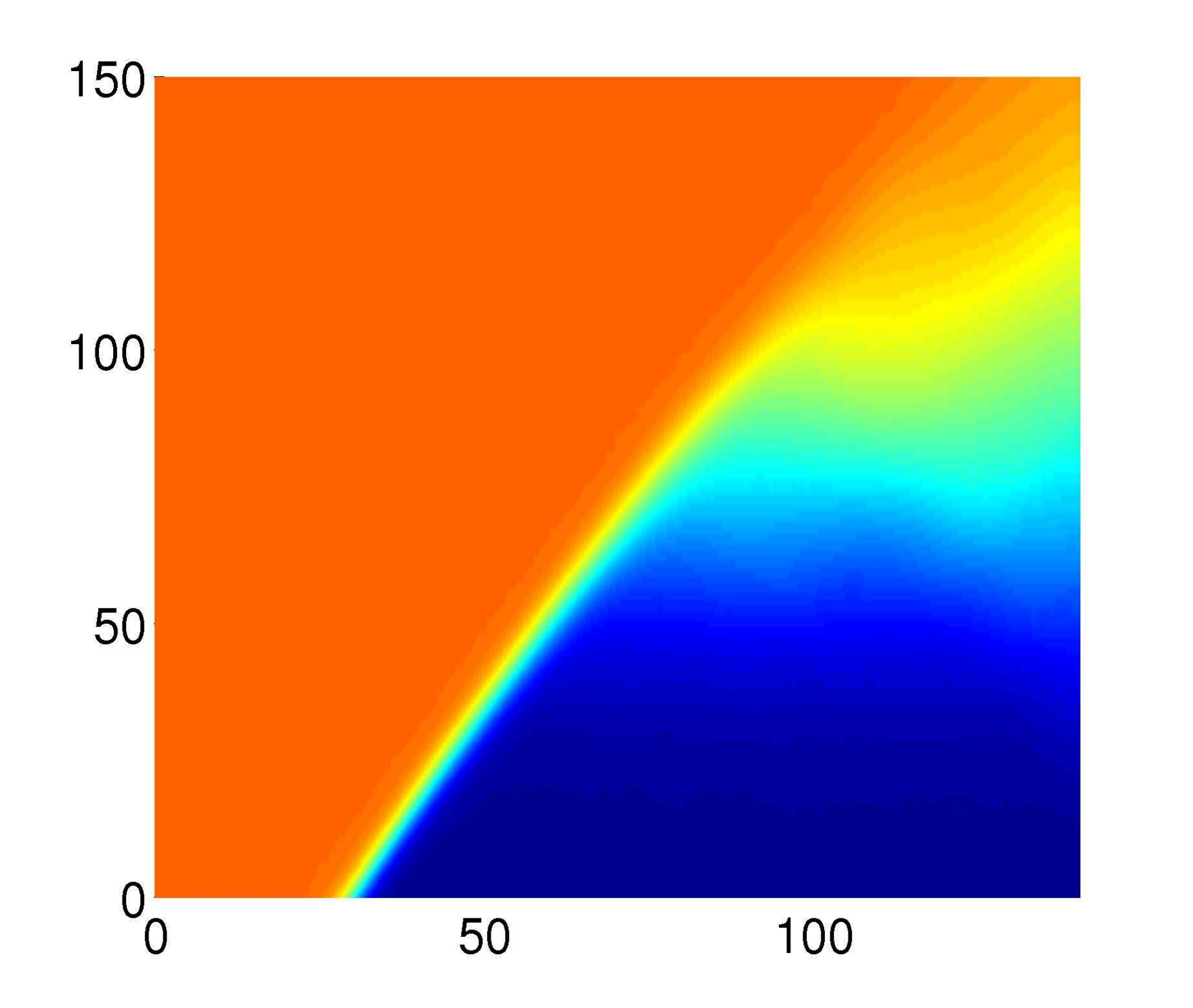}
  \includegraphics*[width=0.45\textwidth]{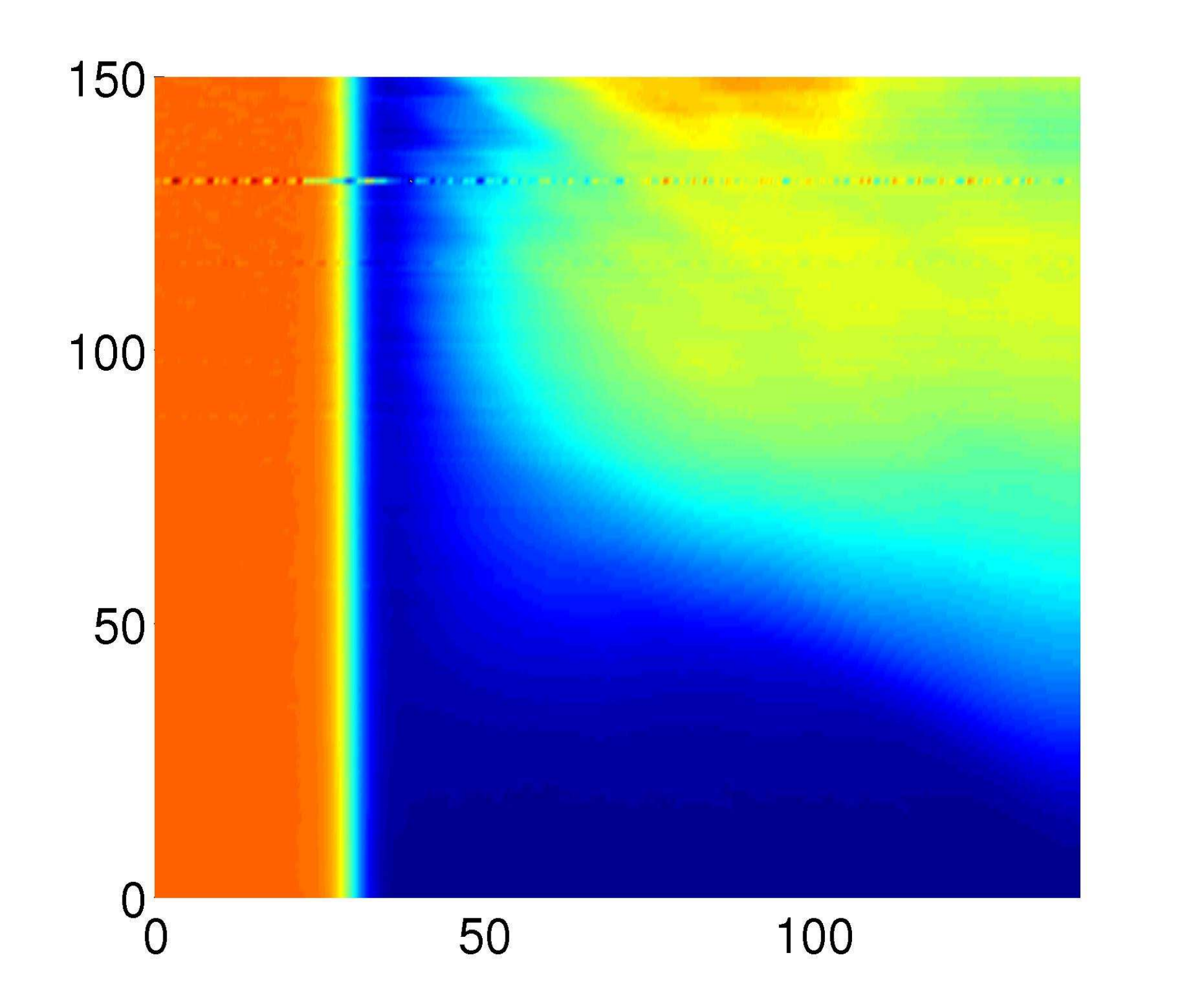}
\end{center}
\caption{Nucleation of travelling waves and annihilation for the
  Nagumo equation with $\alpha=0.1$. In (a) the space-time plot shows
  computations of the \SPDE (no freezing) and in (b) the \SPDAE where
  the wave is frozen. In (c) and (d) are plotted means over
  realizations for the frozen and travelling cases.} 
\label{fig:nuc}
\end{figure}

We now change the parameter $\alpha$ in the nonlinearity 
to $\alpha=0.1$ and illustrate how the \SPDAE approach deals with 
nucleation and extinction of waves. In \figref{fig:nuc} we have plotted 
in (a) a single realization of the \SPDE (so not frozen) showing
nucleation and subsequent extinction ($t\approx 98$) of a
travelling wave. In (b) is plotted a single realization from computing
using the \SPDAE approach. We see the wave is fixed in the domain and
at $t\approx 50$ a wave is nucleated at $x\approx 100$ by the
additive noise. The computations are based on the original wave which
remains fixed until it interacts with the nucleated wave and is
annihilated at $t\approx 94$ when the computations stop when the wave
cease to exist.
In (c) and in (d) we have plotted mean profiles for the \SPDE and the
frozen \SPDAE systems. In each case we see a well defined front from
the averaging and individual nucleations and annihilations are no longer
distinguishable (although in (d) a large solution pollutes the data at
$t\approx 130$).

\section{Discussion}
\label{sec:conc}

We have examined level set based methods and minimization to a
reference function methods to calculate the wave speed of a stochastic
travelling wave. Our numerical results illustrate these give
comparable results. Numerically we saw that for reference functions
with support much smaller than the support of the travelling wave that 
the minimization may fail. Using the minimization technique for the
SPDE (when it is not frozen) is more computationally expensive than
the level set based methods as it requires interpolation at each time
step.

The algorithm described for freezing the wave and solving the SPDAE
has several numerical advantages over simply solving the \SPDE if the
numerical instability issues could be over come. The frozen wave does
not require a large computational domain for long time simulations and
the generation of the noise path is not so computationally expensive.
The cost of of the minimization when the wave is fixed is minimal as
we simply need to compute two inner-products.
However the advection term is nontrivial - and the loss
of numerical stability is a real issue where some realizations fail to
exist as ignoring results where there is numerical blow up may bias
the statistics.

Our investigation of the Nagumo equation has revealed  interesting
and new computational observations that we have not seen reported in
the literature. Although it was known that for Stratonovich noise
increasing noise intensity increases wave speed we have also seen it
increases the support of the wave. In addition increasing the spatial 
correlation decreases the wave speed and decreases the support of the wave.
The reverse is observed for \Ito noise: the noise intensity seems to decrease
the wave speed and correlation length has little influence on the
speed decreases the support of the wave.

For additive noise in the Nagumo equation we see that the wave speed is
increased with the noise intensity like in the multiplicative case --
this is probably because of the small perturbations ahead of the front
that make the wave faster.

\bibliographystyle{siam}
\bibliography{../stw}
\end{document}